\documentclass[reqno, 11pt]{amsart}
\usepackage{amsmath}	
\usepackage{amssymb}	
\usepackage{amsthm}
\usepackage{amsfonts}
\usepackage{latexsym}
\usepackage{mathrsfs}
\usepackage[all,ps,cmtip]{xy}
\usepackage{tabularx}
\usepackage{pict2e}
\usepackage{graphicx}
\usepackage{xypic}
\usepackage{hyperref}
\usepackage{tikz}
\usepackage{makecell}
\usepackage{diagbox}

\usepackage[text={158mm, 235mm},centering]{geometry}
\input diagxy
\input xy
\theoremstyle{plain}
\newtheorem{thm}{Theorem}[section]
\newtheorem{lem}[thm]{Lemma}
\newtheorem{prop}[thm]{Proposition}
\newtheorem{cor}[thm]{Corollary}

\newtheorem*{quest}{Question}

\theoremstyle{definition}
\newtheorem{defn}[thm]{Definition}
\newtheorem{rem}[thm]{Remark}
\newtheorem{example}[thm]{Example}

\DeclareMathOperator{\A}{\mathcal{A}}
\DeclareMathOperator{\CO}{\mathcal{O}}
\DeclareMathOperator{\K}{\mathcal{K}}

\allowdisplaybreaks

\numberwithin{equation}{section}

\begin{document}

\title{Holomorphic Koszul--Brylinski homologies of Poisson blow-ups}

\author[X. Chen]{Xiaojun Chen}
\address{School of Mathematics, Sichuan University, Chengdu 610064, P. R. China}%
\email{xjchen@scu.edu.cn}%

\author[Y. Chen]{Youming Chen }
\address{School of Science, Chongqing University of Technology, Chongqing 400054, P.R. China}%
\email{youmingchen@cqut.edu.cn}%

\author[S. Yang]{Song Yang}
\address{Center for Applied Mathematics, Tianjin University, Tianjin 300072, P.R. China}
\email{syangmath@tju.edu.cn}%

\author[X. Yang]{Xiangdong Yang}
\address{Department of Mathematics, Lanzhou University, Lanzhou 730000, P.R. China}
\email{yangxd@lzu.edu.cn}%


\begin{abstract}
We derive a blow-up formula for holomorphic Koszul--Brylinski
homologies of compact holomorphic Poisson manifolds.
As applications, we investigate the invariance of the $E_{1}$-degeneracy
of the Dolbeault--Koszul--Brylinski spectral sequence under Poisson blow-ups,
and compute the holomorphic Koszul--Brylinski homology for del Pezzo surfaces
and two complex nilmanifolds with holomorphic Poisson structures.
\end{abstract}

\maketitle

\setcounter{tocdepth}{1}
\tableofcontents


\section{Introduction}

Historically, Poisson structures arise from classical mechanics.
In mathematics, the Poisson structures emerge from many fields,
such as generalized complex geometry, geometric representation theory, integrable systems,
and algebraic geometry.
In many situations, the Poisson structures are actually
holomorphic; see \cite{BZ99,Hit06,L-GSX08,Go10,Hit12,CFP16,CGP15} etc..
The study of Poisson structures from the viewpoint of algebraic geometry
can be traced back at least to Bondal \cite{Bon93} and Polishchuk \cite{Pol97}.
More generally, we refer the readers to \cite{Pym18} for an introduction to the
algebraic geometry of holomorphic Poisson structures.
The purpose of this paper is to study holomorphic Poisson structures from
an algebro-geometric point of view, and we focus on the homological aspects
of compact holomorphic Poisson manifolds.

Let $(X,\mathcal{O}_{X})$ be a complex manifold or a scheme of finite type over
the complex number field $\mathbb{C}$.
By a Poisson structure on $X$, we mean a $\mathbb{C}$-bilinear sheaf morphism:
$$
\{-,-\}:\mathcal{O}_{X}\times\mathcal{O}_{X}\longrightarrow \mathcal{O}_{X}
$$
which satisfies the usual axioms for a Poisson bracket, and we call $(X,\{-,-\})$
a holomorphic Poisson manifold or a Poisson scheme.
In particular, holomorphic Poisson structures are closely related to generalized complex geometry.
On the one hand, a holomorphic Poisson structure naturally defines a generalized
complex structure of special type; see \cite{Gu11}.
On the other hand, by Bailey's local classification theorem of generalized complex structures,
each generalized complex manifold is locally equivalent to the product of a symplectic manifold
and a holomorphic Poisson manifold; see \cite{Ba13}.
We refer the readers to \cite{Hit03,Hit06,Fu05,L-GSX08,Go10,CSX10,Gu11,BX15,BCV19}
and references therein for more results on the applications of holomorphic
Poisson structures
in generalized complex geometry and the relationships with other geometries.

Assume that $(X,\{-,-\})$ is a compact holomorphic Poisson manifold of complex dimension $n$,
and let $\pi\in H^{0}(X,\wedge^{2}\mathcal{T}_{X})$ be the holomorphic Poisson bi-vector field
determined by the given Poisson bracket.
For a cohomological study of the Poisson structure of $X$, we have the \emph{holomorphic
Lichnerowicz--Poisson cohomology} $H^{\bullet}(X,\pi)$ defined to be the hypercohomology of
the sheaf complex
$$
\xymatrix{
0\ar[r]&\mathcal{O}_{X} \ar[r]^{b_{\pi}} &\mathcal{T}_{X} \ar[r]^{b_{\pi}} &
\wedge^{2}\mathcal{T}_{X} \ar[r]^{b_{\pi}} & \wedge^{3}\mathcal{T}_{X} \ar[r]^{b_{\pi}} &
\cdots \ar[r]^{b_{\pi}} & \wedge^{n}\mathcal{T}_{X} \ar[r]^{} & 0},
$$
where the differential operator $b_{\pi}(-)=[\pi,-]_S$ is the adjoint
action of $\pi$ with respect to the Schouten bracket; see \cite{Lic77,L-GSX08}.
This cohomology has been widely studied; see, for example,
\cite{HX11,FM12,CGP15,CFP16,PS17,PS19,Hon19} and references therein.
Dually, from a homological point of view, we have the so called {\it holomorphic Koszul--Brylinski complex}:
$$
\xymatrix{
0\ar[r]&\Omega_{X}^{n} \ar[r]^{\partial_{\pi}} &\Omega_{X}^{n-1} \ar[r]^{\partial_{\pi}}
& \Omega_{X}^{n-2} \ar[r]^{\partial_{\pi}} & \Omega_{X}^{n-3}
\ar[r]^{\partial_{\pi}} & \cdots \ar[r]^{\partial_{\pi}}
& \mathcal{O}_{X} \ar[r]^{} & 0},
$$
where $\partial_{\pi}=[\iota_{\pi},\partial]$.
The hypercohomology of the sheaf complex above,
denoted by $H_{\bullet}(X,\pi)$, is called the {\it holomorphic Koszul--Brylinski homology} of $X$.
Most notably, there exists a holomorphic version of Evens--Lu--Weinstein duality
for $H_{\bullet}(X,\pi)$,
which is  a generalization of Serre duality for Dolbeault cohomology; see \cite[Theorem 4.4]{Sti11}.
Furthermore, there is a canonical Fr\"{o}hlicher-type spectral sequence, called the
Dolbeault--Koszul--Brylinski spectral sequence (see Definition \ref{Dol-Poi-spectral-seq}),
which converges to $H_{\bullet}(X,\pi)$.
However, it is not so easy to compute the holomorphic
Koszul--Brylinski homology for a specific holomorphic Poisson manifold.

In algebraic and complex geometry, the blow-up transformation plays a
central role in the  study of algebraic varieties and complex manifolds.
In the Poisson category, it was Polishchuk \cite{Pol97} who first gave the
construction of blow-ups for Poisson schemes.
Polishchuk's construction of blow-up transformations for Poisson schemes
adapts to holomorphic Poisson manifolds without any essential changes.
Our starting point is to understand the homological aspect of holomorphic Poisson manifolds under
a Poisson blow-up transformation.
Particularly, if the holomorphic Poisson structure $\pi$ is trivial,
then the holomorphic Koszul--Brylinski homology is isomorphic to the
Hochschild homology of the complex manifold $X$:
$$
H_{k}(X, 0)
\cong
\bigoplus_{p-q=n-k}H^{q}(X, \Omega_{X}^{p})
\cong
\mathrm{HH}_{n-k}(X).
$$
The blow-up formula for the Hochschild homology has been established in \cite{RYY19}.
To be more specific, suppose $Z \subset X$ is a closed complex manifold of
codimension $c\geq 2$ and $\tilde{X}$ is the blow-up of $X$ along $Z$,
then there exists an isomorphism of Hochschild homologies
$$
 \mathrm{HH}_{n-k}(\tilde{X})\cong  \mathrm{HH}_{n-k}(X) \oplus  \mathrm{HH}_{n-k}(Z)^{\oplus c-1}.
$$
So a natural question that arises now is:

\begin{quest}\label{pro}
For a non-trivial holomorphic Poisson structure,
can we describe explicitly the variance of the holomorphic Koszul--Brylinski homology under a Poisson blow-up?
\end{quest}

Using a sheaf-theoretic approach,
we establish a blow-up formula for holomorphic Koszul--Brylinski homology as follows.

\begin{thm}\label{main-thm}
Suppose $(X, \pi)$ is a compact holomorphic Poisson manifold of complex dimension $n\geq 2$,
and $(Z, \pi|_{Z})\subset (X, \pi)$ is a closed holomorphic Poisson
submanifold of codimension $c\geq 2$ with trivial transverse Poisson structure.
Let $\varphi: \tilde{X}\rightarrow X$ be the blow-up of $X$ along $Z$
with exceptional divisor $E:=\varphi^{-1}(Z)$,
and let $\tilde{\pi}$ be the unique holomorphic Poisson structure on $\tilde{X}$
such that $\varphi$ is a Poisson morphism, i.e., $\varphi_{\star}\tilde{\pi}=\pi$.
Then there exists an isomorphism of holomorphic Koszul--Brylinski homologies:
\begin{equation*}
H_{k}(\tilde{X}, \tilde{\pi})\cong
H_{k}(X, \pi)\oplus \Big( H_{k-1}(E, \tilde{\pi}|_{E})/ \rho^{\star}H_{k-c}(Z, \pi|_{Z})\Big)
\end{equation*}
for any $0\leq k\leq 2n$.
Furthermore, if $Z$ satisfies the $\partial\bar{\partial}$-lemma, then we get
\begin{equation*}
H_{k}(\tilde{X}, \tilde{\pi})
\cong
H_{k}(X, \pi)\oplus H_{k-c}(Z, \pi|_{Z})^{\oplus c-1}.
\end{equation*}
In particular, there exists an isomorphism
$$
H_{k}(\tilde{X}, \tilde{\pi})
\cong
H_{k}(X, \pi)
$$
for $0\leq k\leq c-1$ or $2n-c+1\leq k\leq 2n$.
\end{thm}

Observe that the first page of the Dolbeault--Koszul--Brylinski spectral
sequence of the holomorphic Poisson manifold $(X,\pi)$ is the Dolbeault cohomology:
$$
E_1^{s,t}=H^{t}(X, \Omega_{X}^{n-s})\cong H^{n-s,t}_{\bar{\partial}}(X)
\Longrightarrow
H_{n-s+t}(X,\pi).
$$
The study of the degeneracy of the Dolbeault--Koszul--Brylinski spectral sequence at
$E_{1}$-page may be of independent interest.
As an application of Theorem \ref{main-thm},
we investigate the invariance of such degeneracy under Poisson blow-ups.

\begin{thm}\label{main-thm-2}
With the assumption of Theorem \ref{main-thm}, if $Z$ satisfies the $\partial\bar{\partial}$-lemma then the Dolbeault--Koszul--Brylinski
spectral sequence for $(\tilde{X},\tilde{\pi})$ degenerates at $E_{1}$-page if and only if
it does so for $(X,\pi)$ and $(Z,\pi|_{Z})$.
\end{thm}

It is worth noting that if $X$ is a projective manifold or K\"{a}hler manifold then the closed complex submanifold $Z$ automatically satisfies the $\partial\bar{\partial}$-lemma, and therefore both Theorem \ref{main-thm} and Theorem \ref{main-thm-2} are applicable to these situations.

This paper is organized as follows.
In \S\,\ref{prelim},
we review some basics on holomorphic Poisson manifolds and the holomorphic Koszul--Brylinski homology.
We devote \S\,\ref{Poisson blow-ups and modifications} to Poisson blow-ups and modifications.
In \S\,\ref{Comparison under Poisson projective bundles} we derive the Poisson projective bundles formula for
holomorphic Koszul--Brylinski homology,
a key part of the proof of the main theorems.
In \S\,\ref{Comparison under Poisson blow-ups},
the proofs of the main theorems are given.
In \S\,\ref{examples}, the holomorphic Koszul--Brylinski homologies of
some compact holomorphic Poisson manifolds are computed.
Finally, the Appendix \ref{appendix} gives the Hodge diamond of a six-dimensional complex
nilmanifold in \S\,\ref{a-six-nilfold}.


\subsection*{Acknowledgements}
The last three authors would like to thank the School of Mathematics of Sichuan University
and Tianyuan Mathematical Center in Southwest China for the hospitalities during the winter of 2021.
The second author is very grateful to Professor Zhuo Chen for many useful discussions.
This work is partially supported by the National Nature Science Foundation of China
(Nos. 11890660, 11890663, 12126309, 12126354, 12171351),
the Scientific and Technological Research Program of Chongqing Municipal Education Commission
(Grant No. KJQN202201108),
and the Natural Science Foundation of Tianjin (No. 20JCQNJC02000).


\section{Preliminaries}\label{prelim}

In this section, we review some basic facts on holomorphic Poisson manifolds and
the Koszul--Brylinski homology of holomorphic Poisson manifolds.

\subsection{Holomorphic Poisson manifolds}

Let $X$ be a complex manifold
and let $\CO_X$ be its structure sheaf (i.e., the sheaf of holomorphic functions),
$\Omega_{X}^{p}$ be the sheaf of holomorphic $p$-forms,
$\mathcal{T}_{X}$ be the sheaf of holomorphic vector fields.

\begin{defn}
A complex manifold $X$ is called a holomorphic Poisson manifold if $X$
admits a holomorphic bi-vector field
$\pi\in H^{0}(X, \wedge^{2}\mathcal{T}_{X})$ such that $[\pi,\pi]_{S}=0$,
where $[-,-]_{S}$ is the Schouten bracket.
\end{defn}

Such a holomorphic bi-vector field $\pi$ is called the {\it holomorphic
Poisson bi-vector field} of the holomorphic Poisson manifold $X$,
and the holomorphic Poisson manifold $X$ is also denoted by $(X,\pi)$.
In particular,
for any open subset $U\subset X$,
the ring $\CO_{X}(U)$ is equipped with a Poisson bracket $\{-,-\}$ via $\pi$ such that
for any open subset $V\subset U$ of $X$,
the restriction map
$
\CO_{X}(U)\to \CO_{X}(V)
$
is a morphism of Poisson algebras;
the holomorphic Poisson bi-vector field $\pi$ induces
a sheaf morphism
$\pi^\sharp:\Omega_{X}^{1}\rightarrow\mathcal{T}_{X}$ by contraction
with $\pi$.
For any fixed point $p\in X$, $\mathrm{Rank}(\pi)|_p$ is defined
to be the rank of the linear map $\pi^\sharp|_p$.
Naturally, $\mathrm{Rank}(\pi)|_p$ is even and the following theorem
describes the local structure of a holomorphic Poisson structure (c.f. \cite[Theorem 1.25]{L-GPV13}).

\begin{thm}[Weinstein's splitting theorem]\label{WST}
Let $(X,\pi)$ be a holomorphic Poisson manifold and $p$ is an arbitrary point of $X$.
Suppose $\mathrm{Rank}(\pi)|_p=2r$.
Then there exists a neighborhood $U$ of $p$ with holomorphic coordinates
$\{z_1, \cdots, z_{s}, z_{s+1}, \cdots, z_{s+2r}\}$ centered at $p$, such that  on $U$,
$$
\pi=\sum\limits_{1\leq i,j\leq s}\phi_{ij}(z_1,\cdots,z_s)
  \frac{\partial}{\partial z_i}\wedge \frac{\partial}{\partial z_j}
 +\sum\limits_{i=1}^r\frac{\partial}{\partial z_{s+i}}\wedge \frac{\partial}{\partial z_{s+r+i}}
 $$
where the functions $\phi_{ij}$ are holomorphic functions of $(z_1, \cdots, z_s)$
satisfying  $\phi_{ij}(p)=0$.
Such a local coordinate $\{z_1, \cdots, z_{s}, z_{s+1}, \cdots, z_{s+2r}\}$
is called a splitting coordinate centered at $p$.
\end{thm}

A holomorphic map $f: Y\rightarrow X$ of holomorphic
Poisson manifolds $(Y,\pi_Y)$ and $(X,\pi_X)$ is a Poisson morphism
if and only if $f_{\star}(\pi_{Y}|_{p})=\pi_{X}|_{f(p)}$ for every $p\in Y$;
in this case, we write $f_{\star}\pi_{Y}=\pi_{X}$.
In particular,
let
$\jmath: Z\hookrightarrow X$ be a closed complex submanifold of holomorphic Poisson manifold $X$.
Suppose that $Z$ is also a holomorphic Poisson,
then we say that $Z$ is a {\it closed holomophic Poisson submanifold} of $X$
if the inclusion $\jmath$ is a Poisson morphism.

Analogous to the real case,
there are some intrinsic restrictions on the existence of holomorphic Poisson
submanifolds in a holomorphic Poisson manifold.
For example, due to Weinstein's splitting theorem,
one can prove that each holomorphic symplectic leaf(which is hyper-K\"{a}hler) is a Poisson submanifold,
and every holomorphic Poisson submanifold is the union of some symplectic leaves.
Therefore, if the holomorphic Poisson bi-vector field of $X$
is induced by a holomorphic symplectic form, then
only open subsets of $X$ are Poisson submanifolds.
More precisely, consider a closed holomorphic Poisson
submanifold $\jmath: (Z,\pi|_{Z})\hookrightarrow(X,\pi)$, for any $p\in Z$,
we can choose a neighborhood $U$ of $p$ in $X$ with splitting coordinates
$\{z_1,\cdots,z_s,z_{s+1},\cdots,z_{s+2r}\}$ centered at $p$ satisfying
$$
\pi|_{U}
=\sum\limits_{1\leq i,j\leq s}\phi_{ij}(z_1,\cdots,z_s)
  \frac{\partial}{\partial z_i}\wedge \frac{\partial}{\partial z_j}
   +\sum\limits_{i=1}^r\frac{\partial}{\partial z_{s+i}}\wedge \frac{\partial}{\partial z_{s+r+i}},
$$
such that there exists a neighborhood
$V=U\cap Z=\{z_{1}=0,\cdots,z_c=0\}\subset U$ of $p$ in $Z$
satisfying
$$
(\pi|_{Z})|_{V}
=\sum\limits_{c+1\leq i,j\leq s}\jmath^{\ast}\phi_{ij}
 \frac{\partial}{\partial z_i}\wedge \frac{\partial}{\partial z_j}
 +\sum\limits_{i=1}^r\frac{\partial}{\partial z_{s+i}}\wedge \frac{\partial}{\partial z_{s+r+i}}.
$$

\subsection{Koszul--Brylinski homology}

Koszul--Brylinski homology is introduced independently by
Koszul \cite{Kos84} and Brylinski \cite{Bry88}.
Let $(X, \pi)$ be a holomorphic Poisson manifold.
The Koszul--Brylinski operator of $(X, \pi)$
on the sheaves of holomorphic forms is given as follows:
$$
\partial_{\pi}:= [\iota_{\pi}, \partial]:\Omega_{X}^{p}\longrightarrow \Omega_{X}^{p-1},
$$
where $\Omega_{X}^{p}$ is the sheaf of holomorphic $p$-forms,
$\partial$ is the Dolbeault operator and $\iota_{\pi}$ is the contraction operator
with respect to holomorphic Poisson bi-vector field $\pi$.
\begin{lem}\label{lem1}
Let $(X,\pi)$ be a holomorphic Poisson manifold.
Then we have $\iota_{\pi}\bar{\partial}=\bar{\partial}\iota_{\pi}$ and
$\bar{\partial}\partial_{\pi}+\partial_{\pi}\bar{\partial}=0$
\end{lem}

\begin{proof}
To prove the first statement, it suffices to verify the assertion on an arbitrary coordinate neighborhood of $X$.
Let $(U;z_{1},\cdots,z_{n})$ be a coordinate neighborhood of $X$.
Locally, the holomorphic Poisson bi-vector field can be expressed
as $\pi=\sum_{i,j} c_{ij}\frac{\partial}{\partial z_i}\wedge\frac{\partial}{\partial z_j},$
where $c_{ij}$ are holomorphic functions on $U$.
By definition, for any smooth $(p,q)$-form
$\alpha=fd z_{k_1}\wedge\cdots\wedge
d z_{k_p}\wedge d \bar{z}_{l_1}\wedge\cdots\wedge d \bar{z}_{l_q}$ on $U$, we have
\begin{eqnarray*}
(\iota_{\pi}\bar{\partial}-\bar{\partial}\iota_{\pi})\alpha
&=& \sum_{i,j}c_{ij}\cdot\iota_{\frac{\partial}{\partial z_i}\wedge\frac{\partial}{\partial z_j}}
\bigl(\sum_s \frac{\partial f}{\partial \bar{z}_s}
d \bar{z}_{s}\wedge d z_{k_1}\wedge\cdots\wedge d z_{k_p}
\wedge d \bar{z}_{l_1}\wedge\cdots\wedge d \bar{z}_{l_q}\bigr)\\
&&-\bar{\partial}\bigl(\sum_{i,j}fc_{ij}\cdot\iota_{\frac{\partial}{\partial z_i}\wedge\frac{\partial}{\partial z_j}}
(d z_{k_1}\wedge\cdots\wedge d z_{k_p})\wedge d \bar{z}_{l_1}\wedge\cdots\wedge d \bar{z}_{l_q}\bigr)\\
&=&(-1)^p\sum_s\sum_{i,j} \frac{\partial f}{\partial \bar{z}_s}c_{ij}\cdot
\iota_{\frac{\partial}{\partial z_i}\wedge\frac{\partial}{\partial z_j}}(dz_{k_1}
\wedge\cdots\wedge d z_{k_p})\wedge d\bar{z}_{s}\wedge d \bar{z}_{l_1}\wedge\cdots\wedge d \bar{z}_{l_q}\\
&&-\sum_s\sum_{i,j} \frac{\partial f}{\partial \bar{z}_s} d \bar{z}_{s}\wedge c_{ij}\cdot
\iota_{\frac{\partial}{\partial z_i}\wedge\frac{\partial}{\partial z_j}}
(d z_{k_1}\wedge\cdots\wedge d z_{k_p})\wedge d \bar{z}_{l_1}\wedge\cdots\wedge d \bar{z}_{l_q})\\
&=&0.
\end{eqnarray*}
Equivalently, we get $\bar{\partial}_{\pi}:=[\iota_{\pi},\bar{\partial}]=\iota_{\pi}\bar{\partial}-\bar{\partial}\iota_{\pi}=0$.

For the second statement, since $[\bar{\partial},\partial_{}]=\bar{\partial}\partial_{}+\partial_{}\bar{\partial}=0$,
by the first statement,
we have that
$\bar{\partial}\partial_{\pi}+\partial_{\pi}\bar{\partial}
=[\bar{\partial},\partial_{\pi}]=[\bar{\partial},[\iota_{\pi}, \partial]]
=[[\bar{\partial}, \iota_{\pi}],\partial]+[\iota_{\pi},[\bar{\partial}, \partial]]=0$,
and the lemma is proved.
\end{proof}
According to the Cartan formulae,
we have $\partial_{\pi}^2=0$, 
and
$$
\partial_{\pi}(\alpha\wedge \beta)
=
(\partial_{\pi}\alpha) \wedge \beta
+(-1)^{k} \alpha\wedge (\partial_{\pi}\beta)+(-1)^{k}[\alpha, \beta]_{\partial_{\pi}}
$$
for any  $\alpha\in \Omega_{X}^{k}$ and $\beta\in \Omega_{X}^{l}$.
Here $[-, -]_{\partial_{\pi}}$ is a graded Lie bracket on $\Omega_{X}^{\bullet}$ obtained by Leibniz rule via
\begin{equation}\label{Lie str on forms}
[\alpha, \beta]_{\partial_{\pi}}
:=L_{\pi^{\sharp}(\alpha)}\beta-L_{\pi^{\sharp}(\beta)}\alpha-\partial(\pi(\alpha,\beta)),
\; \forall \; \alpha,\beta\in \Omega_{X}^{1}.
\end{equation}
The  {\it holomorphic Koszul--Brylinski complex} of $X$ is the sheaf complex:
\begin{equation}\label{KB-sheaf-complex}
0\to \Omega_{X}^{n} \stackrel{\partial_{\pi}}{\to}
     \cdots \stackrel{\partial_{\pi}}{\to}\Omega_{X}^{s+1} \stackrel{\partial_{\pi}}{\to}
        \Omega_{X}^{s} \stackrel{\partial_{\pi}}{\to} \Omega_{X}^{s-1} \stackrel{\partial_{\pi}}{\to}
        \cdots\stackrel{\partial_{\pi}}{\to}\CO_X\to 0.
\end{equation}

\begin{defn}
Let $(X, \pi)$ be a holomorphic Poisson manifold.
The $k$-th {\it holomorphic Koszul--Brylinski homology} of $(X,\pi)$ is defined to be
\begin{equation}
H_{k}(X, \pi):=\mathbb{H}^{k}(X, (\Omega_{X}^{\bullet},\partial_{\pi})),
\end{equation}
the $k$-th hyperchomology of the holomorphic Koszul--Brylinski complex.
\end{defn}

\begin{prop}\label{KB-equal1}
Suppose $(X, \pi)$ is a holomorphic Poisson manifold.
Then its holomorphic Koszul--Brylinski complex admits a fine resolution which is
the total complex of the Koszul--Brylinski double complex
$(\A_{X}^{\bullet,\bullet}, \partial_{\pi},\bar{\partial})$,
where $\A_{X}^{p,q}$ is the sheaf of $(p,q)$-forms on $X$.
In particular, the Koszul--Brylinski homology is isomorphic to the
hypercohomology of the associated total complex.
\end{prop}

\begin{proof}
Since the sheaf complex $\A_{X}^{p,\bullet}$ gives rise to a fine resolution of $\Omega_{X}^{p}$,
the assertion follows from the fact that the Koszul--Brylinski operator
$\partial_{\pi}$ commutes with $\bar{\partial}$; see also \cite[Theorem 5.1]{Sti11}.
\end{proof}

This proposition immediately yields the natural morphism of
Koszul--Brylinski homology under Poisson morphisms.

\begin{cor}
Suppose that $f: (Y,\pi_Y)\rightarrow (X,\pi_X)$ is a Poisson morphism of holomorphic Poisson manifolds.
Then the pullback of differential forms naturally induces
a morphism of the holomorphic Koszul--Brylinski homologies
$$
f^{\star}: H_{k}(X, \pi_{X})\longrightarrow H_{k}(Y, \pi_{Y}).
$$
\end{cor}

\begin{proof}
Note that on the space of $(p,q)$-forms,
we have
$$
f^{\star}\circ \partial_{\pi_{X}}= f^{\star}\circ\partial_{f_{\star}\pi_{Y}}
= \partial_{\pi_{Y}}\circ f^{\star}\;\; \textrm{and}\;
f^{\star}\circ \bar{\partial}=\bar{\partial}\circ f^{\star}.
$$
Hence, the corollary follows immediately from Proposition \ref{KB-equal1}.
\end{proof}

By a result of Sti\'{e}non \cite[Theorem 6.4]{Sti11},
the holomorphic Evens--Lu--Weinstein pairing on the holomorphic Koszul--Brylinski
homology is non-degenerate.
More precisely,
if $(X, \pi)$ is a compact holomorphic Poisson manifold of complex dimension $n$,
then there is an isomorphism
\begin{equation}\label{Serre-Poincare-duality}
H_{2n-k}(X,\pi)\cong H_{k}(X,\pi)
\end{equation}
for $0\leq k\leq 2n$.
In the dual aspect, there exists a {\it holomorphic Lichnerowicz--Poisson complex}
$(\wedge^{\bullet}\mathcal{T}_{X},b_{\pi})$:
\begin{equation*}
0\to \CO_X\stackrel{b_{\pi}}{\to}
     \cdots \stackrel{b_{\pi}}{\to}\wedge^{s-1}\mathcal{T}_{X} \stackrel{b_{\pi}}{\to}
           \wedge^{s}\mathcal{T}_{X} \stackrel{b_{\pi}}{\to} \wedge^{s+1}\mathcal{T}_{X} \stackrel{b_{\pi}}{\to}
                    \cdots\stackrel{b_{\pi}}{\to}\wedge^{n}\mathcal{T}_{X} \to 0
\end{equation*}
where $b_{\pi}(-)=[\pi,-]_S$.
The $k$-th hyperchomology of $(\wedge^{\bullet}\mathcal{T}_{X},b_{\pi})$ is called the
$k$-th \emph{holomorphic Lichnerowicz--Poisson cohomology}, i.e.,
$$
H^{k}(X, \pi):=\mathbb{H}^{k}(X, (\wedge^{\bullet}\mathcal{T}_{X},b_{\pi})).
$$
Assume that $X$ admits a holomorphic volume form $\omega\in \Gamma(X,\Omega^n_X)$.
Then there is the natural morphism of sheaves
$$
\iota_{(-)}\omega:
\wedge^{s}\mathcal{T}_{X}
\longrightarrow
\Omega_{X}^{n-s}
$$
for each $s\in \{0,1, \cdots, n\}$.
However, it does not induce a morphism of sheaf complexes
between
$(\wedge^{\bullet}\mathcal{T}_{X},b_{\pi})$ and $(\Omega_{X}^{\bullet},\partial_{\pi})$.
The reason lies in the fact that the diagram
$$
\xymatrixcolsep{3pc}
\xymatrix{
\wedge^{s}\mathcal{T}_{X} \ar[d]_{b_{\pi}} \ar[r]^{\iota_{(-)}\omega} & \Omega^{n-s}_X\ar[d]^{\partial_\pi}\\
\wedge^{s+1}\mathcal{T}_{X} \ar[r]^{\iota_{(-)}\omega} & \Omega^{n-s-1}_X
}
$$
is not commutative in general.
This motivates the following definition.

\begin{defn}[{c.f. \cite{We97,BZ99}}]
A holomorphic Poisson manifold $(X,\pi)$ is called {\it unimodular}
if there is a holomorphic volume form $\omega$
such that the morphism  $\iota_{(-)}\omega$
induces a morphism of sheaf complexes from
$(\wedge^{\bullet}\mathcal{T}_{X},b_{\pi})$ to $(\Omega^{\bullet}_{X},\partial_{\pi})$.
\end{defn}

Equivalently,
a holomorphic Poisson manifold $(X,\pi)$ is unimodular if and only if $\partial_\pi\omega=0$,
or the modular vector field, introduced by Weinstein \cite{We97} and
Brylinski--Zuckerman \cite{BZ99}, vanishes.
In particular, we have

\begin{prop}[{\cite[Proposition 4.7]{Sti11}}]\label{dual}
If the holomorphic Poisson manifold $(X,\pi)$ is unimodular,
then there is an isomorphism
$$
H_{k}(X, \pi)\cong H^{2n-k}(X, \pi),
$$
for any $k\in \mathbb{Z}$, where $n=\mathrm{dim}_{\mathbb{C}}X$.
\end{prop}


\section{Blow-ups and modifications in the Poisson
category}\label{Poisson blow-ups and modifications}

In this section, we give a rapid review on the blow-ups and modifications
in the holomorphic Poisson category.

\subsection{Poisson blow-ups}

Given a complex manifold $X$ and a closed complex
submanifold $\jmath:Z \hookrightarrow X$ with complex codimension $c\geq 2$.
Let $\varphi: \tilde{X}\rightarrow X$ be the blow-up of $X$ along $Z$.
Then the holomorphic map
$$
\varphi: \tilde{X}-E \to X-Z
$$
is biholomorphic, where $E:=\varphi^{-1}(Z)$
is the exceptional divisor, which is the projective bundle of the normal bundle of $Z$ in $X$.
Moreover, we have a commutative diagram

\begin{equation}\label{blowup_diagram}
\vcenter{
\xymatrixcolsep{3pc}
\xymatrix{
E \ar[d]_{\rho:=\varphi|_E } \ar@{^{(}->}[r]^{\tilde{\jmath}} & \tilde{X}\ar[d]^{\varphi}\\
 Z \ar@{^{(}->}[r]^{\jmath} & X.
}}
\end{equation}

In the Poisson category,
if $X$ is a holomorphic Poisson manifold and $Z$ is a closed holomorphic Poisson submanifold of $X$,
then the existence of the holomorphic
Poisson structure on the complex blow-up $\tilde{X}$ is not an unconditional result.
In fact, there exist some restrictions on the existence of the holomorphic Poisson structure on $\tilde{X}$.
Let us recall the result which was originally studied by Polishchuk \cite{Pol97}.
Assume $(Z,\pi|_{Z})$ is a closed holomorphic Poisson submanifold of $(X,\pi)$.
Then for any point $z\in Z$, the conormal space $N_{z}^{\ast}Z$
is a Lie algebra induced by the bracket $\eqref{Lie str on forms}$,
or equivalently, the normal space $N_{z}Z$ admits a linear Poisson structure
which defines the transverse Poisson structure $\pi_{N}\in\Gamma(Z,N^{\ast}Z\otimes\wedge^{2}NZ)$.

\begin{defn}
The transverse Poisson structure $\pi_{N}$ of a closed holomorphic Poisson submanifold
$(Z,\pi|_{Z})$ in $(X,\pi_{X})$ is said to be \emph{degenerate} if, for any point $z\in Z$, the map
\begin{eqnarray*}
\wedge^{3} N_{z}^{\ast}Z &\longrightarrow& S^2N_{z}^{\ast}Z \\
\alpha\wedge\beta\wedge\gamma&\longmapsto&
[\alpha,\beta]\gamma+[\beta,\gamma]\alpha+[\gamma,\alpha]\beta
\end{eqnarray*}
is identically to zero.
\end{defn}

It follows from \cite[Proposition 8.1]{Pol97} that a degenerate Lie algebra is either abelian
or isomorphic to the Lie algebra $\mathrm{Span}\{e_{1},\cdots, e_{c-1},f\}$
with Lie bracket $[e_{i},e_{j}]=0, [f, e_{i}]=e_{i}$.

\begin{example}[{c.f. \cite[\S\,2.5.2]{Pym18}}]
Let $\pi$ be a holomorphic Poisson bi-vector field on $\mathbb{C}^2$,
and $\mathrm{Bl}_{o}\mathbb{C}^2\stackrel{\varphi}{\to} \mathbb{C}^2$
the blow-up of $\mathbb{C}^2$ at the origin $o=(0,0)\in \mathbb{C}^{2}$.
Choose coordinates $z_1, z_2$, and suppose
$\{z_1,z_2\}=f(z_1,z_2)$ for some holomorphic function $f$.
Set
$$
u=\varphi^\ast(z_1), \,\,\,\mathrm{and}\,\,\, v
=\frac{\varphi^{\ast}(z_2)}{\varphi^{\ast}(z_1)}=\varphi^{\ast}(z_1^{-1}z_2).
$$
Suppose we can define a holomorphic Poisson bracket on
$\mathrm{Bl}_{o}\mathbb{C}^2$
which is compatible with the one determined by $\pi$ on $\mathbb{C}^2$;
then
\begin{eqnarray*}
\{u,v\}
&=& \{\varphi^\ast(z_1), \varphi^\ast(z_1^{-1}z_2)\}
=\varphi^\ast\{z_1,z_1^{-1}z_2\} =\varphi^\ast(z_1^{-1}f(z_1,z_2)) \\
&=& u^{-1}f(u,uv)=u^{-1}(f(0,0)+ug(u,v)),
\end{eqnarray*}
where $g$ is holomorphic near the locus $u=0$.
Therefore the holomorphic Poisson bracket given by $\{z_1,z_2\}=f(z_1,z_2)$ on $\mathbb{C}^2$
can be lifted to $\mathrm{Bl}_{o}\mathbb{C}^2$ if and only if $f(0,0)=0$.
\end{example}

Now, let us return to the construction of Poisson blow-ups.
The blow-up of a Poisson scheme was originally clarified in the work of Polishchuk \cite{Pol97}.
Here, we review the blow-up of holomorphic Poisson manifolds along
closed holomorphic Poisson submanifolds; see also \cite[Section 2]{Fu05}.

\begin{prop}[{\cite[Propositions 8.2 \& 8.3]{Pol97}
or \cite[Proposition 3.15]{BCV19}}]\label{Constr-Poisson-bp}
Let $(X,\pi)$ be a holomorphic Poisson manifold.
Suppose $\jmath:(Z,\pi|Z) \hookrightarrow (X,\pi)$ is a closed holomorphic Poisson submanifold.
If the associated transverse Poisson structure $\pi_{N}$ vanishes,
then the following statements hold:
\begin{itemize}
\item[(i)] there exists a unique holomorphic Poisson structure $\tilde{\pi}$ on
$\tilde{X}$ such that $\varphi$ is Poisson morphism (i.e., $\varphi_{\star} \tilde{\pi}=\pi$);
\item[(ii)] $E$ is a holomorphic Poisson manifold such that $\varphi|_{E}:E\rightarrow Z$ is a Poisson morphism;
\item[(iii)]the diagram \eqref{blowup_diagram} of holomorphic Poisson manifolds is commutative.
\end{itemize}
\end{prop}

\subsection{Poisson modifications}
This subsection is devoted to the study of the behavior of the holomorphic
Koszul--Brylinski homology under Poisson modifications of compact holomorphic Poisson manifolds.
Recall that a {\it modification} of compact complex manifolds
is a holomorphic map $\psi:Y\rightarrow X$ of compact complex manifolds
satisfying:
\begin{itemize}
  \item [(i)]$\dim\, Y=\dim\, X$; and
  \item [(ii)] there is an analytic subset $S\subset X$ of
  codimension $\geq2$ such that the restriction
  $$
  \psi:Y-\psi^{-1}(S)\longrightarrow X-S
  $$
  is biholomorphic.
\end{itemize}

\begin{defn}
A {\it Poisson modification} is a Poisson morphism
$\psi:(Y,\pi_{Y})\rightarrow (X,\pi_{X})$ of compact holomorphic
Poisson manifolds $(X, \pi_{X})$ and $(Y,\pi_{Y})$
such that $\psi$ is also a modification of compact complex manifolds.
\end{defn}

Note that the holomorphic Poisson blow-ups are important examples of Poisson modifications.
To study the behavior of the holomorphic Koszul--Brylinski homology
under Poisson modifications of compact holomorphic Poisson manifolds,
we need to reinterpret the Koszul--Brylinski homology in terms of currents.
Let $(X,\pi)$ be a holomorphic Poisson manifold,
and $\mathcal{C}_{X}^{s,t}$ be the sheaf of $(s,t)$-currents on $X$.
Then the operators $\partial_{\pi}$ and $\bar{\partial}$ naturally induce
the dual operators $\partial_{\pi}^{\star}$ and $\bar{\partial}^{\star}$
acting on $\mathcal{C}_{X}^{s,t}$, respectively.
Since $\partial_{\pi}^{\star}$ commutes with $\bar{\partial}^{\star}$,
we obtain a double complex
$(\mathcal{C}_{X}^{\bullet,\bullet},\partial_{\pi}^{\star},\bar{\partial}^{\star})$.
In particular, there exists a natural morphism of double complexes
\begin{equation}\label{double-map1}
\tau_{X}:
(\mathcal{A}_{X}^{\bullet,\bullet},\partial_{\pi},\bar{\partial})
\hookrightarrow
(\mathcal{C}_{X}^{\bullet,\bullet},\partial_{\pi}^{\star},\bar{\partial}^{\star}).
\end{equation}
Denote by  $H_{k}^{C}(X, \pi)$ the $k$-hypercohomology of the total complex
of the double complex
$(\mathcal{C}_{X}^{\bullet,\bullet},\partial_{\pi}^{\star},\bar{\partial}^{\star})$.

\begin{lem}\label{KB-equal2}
The natural morphism $\tau_{X}$ induces an isomorphism
$$
\tau_{X}:
H_{k}(X, \pi)
\longrightarrow
H_{k}^{C}(X, \pi),
$$
for any $k\in \mathbb{Z}$.
\end{lem}

\begin{proof}
To prove the assertion, it suffices to verify that
$(\mathcal{A}_{X}^{\bullet,\bullet},\partial_{\pi},\bar{\partial})$
is quasi-isomorphic to
$(\mathcal{C}_{X}^{\bullet,\bullet},\partial_{\pi}^{\star},\bar{\partial}^{\star})$
under the morphism $\tau_{X}$.
By the spectral sequence theory for double complexes,
there exists a sequence $\{E_{r},d_{r}\}$ for
$(\mathcal{A}_{X}^{\bullet,\bullet},\partial_{\pi},\bar{\partial})$
such that
$$
E_{1}=H^{\bullet}(\mathcal{A}_{X}^{\bullet,\bullet},\bar{\partial})
=H^{\bullet,\bullet}_{\bar{\partial}}(X)\Longrightarrow E_{\infty}=H_{\bullet}(X,\pi).
$$
Similarly, the double complex
$\{\mathcal{C}_{X}^{\bullet,\bullet},\partial_{\pi}^{\star},\bar{\partial}^{\star}\}$
admits a spectral sequence $(\tilde{E}_{r},\tilde{d}_{r})$ satisfying
$$
\tilde{E}_{1}=H^{\bullet}(\mathcal{C}_{X}^{\bullet,\bullet},\bar{\partial}^{\star})
\Longrightarrow \tilde{E}_{\infty}=H^{C}_{\bullet}(X,\pi).
$$
Observe that \eqref{double-map1} induces a morphism of spectral sequences
$$
\tau_{X,r}:\{E_{r},d_{r}\}\longrightarrow\{\tilde{E}_{r},\tilde{d}_{r}\}.
$$
Since the natural inclusion
$\tau_{X}: (\mathcal{A}_{X}^{p,\bullet},\bar{\partial})\hookrightarrow
(\mathcal{C}_{X}^{p,\bullet},\bar{\partial}^{\star})$ is a quasi-isomorphism,
we get that the induced map $\tau_{X,1}:E_{1}\rightarrow \tilde{E}_{1}$ is an isomorphism
and therefore
$E_{\infty}\cong \tilde{E}_{\infty}$ under $\tau_{X}$.
This implies that \eqref{double-map1} is a quasi-isomorphism and the proof is completed.
\end{proof}

We are ready to present the following comparison theorem for holomorphic
Koszul--Brylinski homology under Poisson modifications.
\begin{thm}\label{KB-injective}
Let $f: (Y,\pi_{Y})\rightarrow (X,\pi_{X})$ be a Poisson modification
of compact holomorphic Poisson manifolds.
Then the natural morphism
$$
f^{\star}: H_{k}(X, \pi_{X})\longrightarrow H_{k}(Y, \pi_{Y})
$$
is injective, for any $k\in \mathbb{Z}$.
\end{thm}

\begin{proof}
Since $f$ is a Poisson morphism, by definition, we have $f_{\star}\pi_{Y}=\pi_{X}$.
This implies
\begin{equation}\label{op-com}
f^{\star}\circ \partial_{\pi_{X}}= \partial_{\pi_{Y}}\circ f^{\star}\,\,\,
\textup{and}\,\,\, f_{\star}\circ \partial_{\pi_{Y}}^{\star}= \partial_{\pi_{X}}^{\star}\circ f_{\star}.
\end{equation}
In particular, we obtain a diagram
\begin{equation}\label{com-dia-cur}
\vcenter{
\xymatrixcolsep{3pc}
\xymatrix{
  (\Gamma(X,\mathcal{A}_{X}^{\bullet,\bullet}),\partial_{\pi_{X}},
  \bar{\partial}) \ar[d]_{f^{\star}} \ar@{^{(}->}[r]^{\tau_{X}}_{}
  & (\Gamma(X,\mathcal{C}_{X}^{\bullet,\bullet}),\partial_{\pi_{X}}^{\star},\bar{\partial}^{\star})   \\
  (\Gamma(Y,\mathcal{A}_{Y}^{\bullet,\bullet}),\partial_{\pi_{Y}},\bar{\partial})
  \ar@{^{(}->}[r]^{\tau_{Y}}_{} &
  (\Gamma(Y,\mathcal{C}_{Y}^{\bullet,\bullet}),\partial_{\pi_{Y}}^{\star},\bar{\partial}^{\star})
  \ar[u]^{f_{\star}} .}}
\end{equation}
However, it is not \emph{a priori} clear that the diagram \eqref{com-dia-cur} is commutative.
We now show the commutativity of \eqref{com-dia-cur}.
As $f$ is a modification of compact complex manifolds,
its degree is $1$; moreover, $f$ is a biholomorphism outside of two sets with
Lebesgue measure zero.
As a result, let $\alpha$ be a differential $k$-form on $X$, then we have
\begin{eqnarray*}
\langle f_{\star}\circ \tau_{Y}\circ f^{\star}(\alpha),\beta \rangle
&=&
\int_{X} (f_{\star}\circ \tau_{Y}\circ f^{\star}(\alpha)) \wedge \beta \\
&=&\int_{Y} f^{\star}(\alpha\wedge \beta)=\int_{X} \alpha\wedge \beta \\
&=& \langle \tau_{X}(\alpha),\beta \rangle,
\end{eqnarray*}
where $\beta$ is an arbitrary differential $(2n-k)$-form on $X$.
It follows that $\tau_{X}(\alpha)=f_{\star}\circ \tau_{Y} \circ f^{\star}(\alpha)$;
see the proof of \cite[Theorem 12.9]{Dem12}.
Combining it with \eqref{op-com} yields that \eqref{com-dia-cur} is a commutative diagram.
Applying Lemma \ref{KB-equal2} to $X$ and $Y$,
we obtain two natural isomorphisms
$$
\tau_{Y}:
H_{k}(Y, \pi_{Y})
\longrightarrow
H_{k}^{C}(Y, \pi_{Y})
\,\,\,\textup{and}\,\,\,
\tau_{X}:
H_{k}(X, \pi_{X})
\longrightarrow
H_{k}^{C}(X, \pi_{X}).
$$
Consequently, we obtain a commutative diagram
\begin{equation*}
\xymatrixcolsep{3pc}
\xymatrix{
  H_{k}(X, \pi_{X}) \ar[d]_{f^{\star}} \ar[r]^{\tau_{X}}_{\simeq} & H_{k}^{C}(X, \pi_{X})   \\
  H_{k}(Y, \pi_{Y})
  \ar[r]^{\tau_{Y}}_{\simeq} & H_{k}^{C}(Y, \pi_{Y}) \ar[u]^{f_{\star}}.}
\end{equation*}
and hence the morphism
$$
f^{\star}: H_{k}(X, \pi_{X})
\longrightarrow H_{k}(Y, \pi_{Y})
$$
is injective.
\end{proof}

\section{Comparison under Poisson projective bundles}\label{Comparison under Poisson projective bundles}
The purpose of this section is to establish the following projective bundle formula for holomorphic
Koszul--Brylinski homology.
\begin{thm}\label{Poisson-pjbd-formula}
Suppose $(Z,\pi)$ is a compact holomorphic Poisson manifold.
Let $\rho: E\rightarrow Z$ be the projective bundle of a holomorphic vector bundle of rank $c\geq 2$ on $Z$.
If $Z$ satisfies the $\partial\bar{\partial}$-lemma and $\tilde{\pi}$ is a holomorphic Poisson structure on $E$ such that $\rho_{\star}\tilde{\pi}=\pi$,
then there is an isomorphism of Koszul--Brylinski homology as $\mathbb{C}$-vector spaces:
$$
H_{k+1-c}(Z, \pi)^{\oplus\; c}\cong H_{k}(E, \tilde{\pi}),
$$
for any $k\in \mathbb{Z}$.
\end{thm}

To illustrate the basic idea of the proof of the theorem above,
we consider the case of $\dim_{\mathbb{C}}\,Z=2$ and $c=3$.
Consider the first Chern class of the tautological line bundle over $E$:
$$\textbf{h}=c_{1}(\mathcal{O}_{E}(1))\in H^{1,1}_{\bar{\partial}}(E).$$
Set $A_{Z}^{s,t}:=\Gamma(Z, \mathcal{A}_{Z}^{s,t})$ be the space of differential $(s,t)$-forms.
Observe that $H_{k}(E,\tilde{\pi})$ is equal to the $k$-th
total cohomology of the double complex
$G=(A_{E}^{\bullet, \bullet},\partial_{\tilde{\pi}},\bar{\partial})$, whereas
$H_{k-2}(Z, \pi)^{\oplus\; 3}$ is the $k$-th total cohomology of the double complex
$$L=\biggl(\bigoplus_{i=0}^{2} A_{Z}^{\bullet, \bullet}[-2+i, -i], \partial_{{\pi}},\bar{\partial}\biggr).$$
According to the standard spectral sequence theory for double complexes,
we have a spectral sequence $\{(\mathcal{G}^{\bullet,\bullet}_{r},d_{r})\}$ associated to $G$ such that
$$
\mathcal{G}^{\bullet,\bullet}_{1}=H^{\bullet,\bullet}_{\bar{\partial}}(E)
\Longrightarrow H_{\bullet}(E,\tilde{\pi}).
$$
Similarly, for the double complex $L$,
there exists a spectral sequence $\{(\mathcal{L}^{\bullet,\bullet}_{r},\bar{d}_{r})\}$ satisfying
$$
\mathcal{L}^{\bullet,\bullet}_{1}=\bigoplus_{i=0}^{2} H_{\bar{\partial}}^{\bullet, \bullet}(Z)[-2+i, -i]
\Longrightarrow
H_{\bullet-2}(Z,\pi)^{\oplus 3}.
$$
Note that there exists a well-defined map of bi-graded $\mathbb{C}$-vector spaces (see the figure below):
$$
\Psi:=\sum_{i=0}^{2}h^{i}\wedge \rho^{\star}(-):
\bigoplus_{i=0}^{2} A_{Z}^{\bullet, \bullet}[-2+i, -i]
\longrightarrow  A_{E}^{\bullet, \bullet}.
$$
Since $\partial_{\tilde{\pi}}$ is not a derivation,
it does not commute with $\Psi$ and therefore
$\Psi$ can not give rise to a morphism between the double complexes $G$ and $L$.
Recall that a compact complex manifold $X$ satisfies the \emph{$\partial\bar{\partial}$-lemma}, if the equation
$$\ker\,\partial\cap\ker\,\bar{\partial}\cap\mathrm{im}\,d=\mathrm{im}\,\partial\bar{\partial}$$
holds for the double complex $(A^{\bullet,\bullet}_{X}, \partial, \bar{\partial})$ (cf. \cite{DGMS75}).
Under the assumption that $Z$ satisfies the
$\partial\bar{\partial}$-lemma, it is noteworthy that $\Psi$ induces a morphism
$\Psi_{1}:\mathcal{G}^{\bullet,\bullet}_{1}\rightarrow\mathcal{L}^{\bullet,\bullet}_{1}$
which commutes with the differentials $d_{1}$ and $\bar{d}_{1}$.
Consequently, we get a well-defined morphism of spectral sequences
$$\Psi_{r}:(\mathcal{G}^{\bullet,\bullet}_{r},d_{r})
\longrightarrow(\mathcal{L}^{\bullet,\bullet}_{r},\bar{d}_{r}).$$
In particular, by the projective bundle formula of Dolbeault cohomology,
we conclude that $\Psi_{1}$ is an isomorphism, and so is the $\Psi_{\infty}$.
\begin{center}
 \begin{tikzpicture}[scale=0.3]
\draw  (-19,0) node[left] {$0$}
         (-19,4) node[left] {$0$}
            (-19,8) node[left] {$0$};
\draw  (-16.2,-3.3) node[left] {$0$}
           (-12.2,-3.3) node[left] {$0$}
                      (-8.2,-3.3) node[left] {$0$};
\draw[-latex] (-17,-2.5)-- (-17,-1);
\draw[-latex] (-13,-2.5)-- (-13,-1);
\draw[-latex] (-9,-2.5)-- (-9,-1);
\draw[-latex] (-17,1.5)-- (-17,3);
\draw[-latex] (-13,1.5)-- (-13,3);
\draw[-latex] (-9,1.5)-- (-9,3);
\draw[-latex] (-17,5.5)-- (-17,7);
\draw[-latex] (-13,5.5)-- (-13,7);
\draw[-latex] (-9,5.5)-- (-9,7);
\draw[-latex] (-17,9.5)-- (-17,11);
\draw[-latex] (-13,9.5)-- (-13,11);
\draw[-latex] (-9,9.5)-- (-9,11);
\draw  (-16.2,11.8) node[left] {$0$}
           (-12.2,11.8) node[left] {$0$}
                      (-8.2,11.8) node[left] {$0$};
  \draw[-latex] (-19.2,0)--(-18,0) node[right] {$A_{Z}^{2,0}$};
  \draw[-latex] (-15.2,0)-- (-14,0) node[right] {$A_{Z}^{2,1}$};
   \draw[-latex] (-11.2,0)-- (-10,0) node[right] {$A_{Z}^{2,2}$};
           \draw[-latex] (-7.2,0)-- (-6,0) node[right] {$0$};
  \draw[-latex] (-19.2,4)--(-18,4) node[right] {$A_{Z}^{1,0}$};
  \draw[-latex] (-15.2,4)-- (-14,4) node[right] {$A_{Z}^{1,1}$};
   \draw[-latex] (-11.2,4)-- (-10,4) node[right] {$A_{Z}^{1,2}$};
      \draw[-latex] (-7.2,4)-- (-6,4) node[right] {$0$};
  \draw[-latex] (-19.2,8)--(-18,8) node[right] {$A_{Z}^{0,0}$};
  \draw[-latex] (-15.2,8)-- (-14,8) node[right] {$A_{Z}^{0,1}$};
   \draw[-latex] (-11.2,8)-- (-10,8) node[right] {$A_{Z}^{0,2}$};
      \draw[-latex] (-7.2,8)-- (-6,8) node[right] {$0$};
  \draw[red] (-19,-1)--(-19,9)-- (-7,9)--(-7,-1)--(-19,-1);
  \draw[-latex,red] (-7,5)-- (1,11);
\draw[-latex] (3,-2.5)-- (3,-1);
\draw[-latex] (7,-2.5)-- (7,-1);
\draw[-latex] (11,-2.5)-- (11,-1);
\draw[-latex] (15,-2.5)-- (15,-1);
\draw[-latex] (19,-2.5)-- (19,-1);
 \draw[-latex] (0.8,0)-- (2,0) node[right] {$A_{E}^{4,0}$};
  \draw[-latex] (4.8,0)-- (6,0) node[right] {$A_{E}^{4,1}$};
   \draw[-latex] (8.8,0)-- (10,0) node[right] {$A_{E}^{4,2}$};
    \draw[-latex] (12.8,0)-- (14,0) node[right] {$A_{E}^{4,3}$};
     \draw[-latex] (16.8,0)-- (18,0) node[right] {$A_{E}^{4,4}$};
           \draw[-latex] (20.8,0)-- (22,0) node[right] {$0$};
\draw[-latex] (3,1)-- (3,3);
\draw[-latex] (7,1)-- (7,3);
\draw[-latex] (11,1)-- (11,3);
\draw[-latex] (15,1)-- (15,3);
\draw[-latex] (19,1)-- (19,3);
  \draw[-latex] (0.8,4)-- (2,4) node[right] {$A_{E}^{3,0}$};
  \draw[-latex] (4.8,4)-- (6,4) node[right] {$A_{E}^{3,1}$};
   \draw[-latex] (8.8,4)-- (10,4) node[right]{$A_{E}^{3,2}$};
    \draw[-latex] (12.8,4)-- (14,4) node[right]{$A_{E}^{3,3}$};
     \draw[-latex] (16.8,4)-- (18,4) node[right]{$A_{E}^{3,4}$};
      \draw[-latex] (20.8,4)-- (22,4) node[right]{$0$};

\draw[-latex] (3,5)-- (3,7);
\draw[-latex] (7,5)-- (7,7);
\draw[-latex] (11,5)-- (11,7);
\draw[-latex] (15,5)-- (15,7);
\draw[-latex] (19,5)-- (19,7);
  \draw[-latex] (0.8,8)-- (2,8) node[right] {$A_{E}^{2,0}$};
  \draw[-latex] (4.8,8)-- (6,8) node[right] {$A_{E}^{2,1}$};
   \draw[-latex] (8.8,8)-- (10,8) node[right] {$A_{E}^{2,2}$};
    \draw[-latex] (12.8,8)-- (14,8) node[right] {$A_{E}^{2,3}$};
     \draw[-latex] (16.8,8)-- (18,8) node[right] {$A_{E}^{2,4}$};
      \draw[-latex] (20.8,8)-- (22,8) node[right] {$0$};
 \draw[-latex] (3,9)-- (3,11);
\draw[-latex] (7,9)-- (7,11);
\draw[-latex] (11,9)-- (11,11);
\draw[-latex] (15,9)-- (15,11);
\draw[-latex] (19,9)-- (19,11);
 \draw[-latex] (0.8,12)-- (2,12) node[right] {$A_{E}^{1,0}$};
  \draw[-latex] (4.8,12)-- (6,12) node[right] {$A_{E}^{1,1}$};
   \draw[-latex] (8.8,12)-- (10,12) node[right] {$A_{E}^{1,2}$};
    \draw[-latex] (12.8,12)-- (14,12) node[right]{$A_{E}^{1,3}$};
     \draw[-latex] (16.8,12)-- (18,12) node[right]{$A_{E}^{1,4}$};
      \draw[-latex] (20.8,12)-- (22,12) node[right]{$0$};
 \draw[-latex] (3,13)-- (3,15);
\draw[-latex] (7,13)-- (7,15);
\draw[-latex] (11,13)-- (11,15);
\draw[-latex] (15,13)-- (15,15);
\draw[-latex] (19,13)-- (19,15);
 \draw[-latex] (0.8,16)-- (2,16) node[right] {$A_{E}^{0,0}$};
  \draw[-latex] (4.8,16)-- (6,16) node[right] {$A_{E}^{0,1}$};
   \draw[-latex] (8.8,16)-- (10,16) node[right] {$A_{E}^{0,2}$};
    \draw[-latex] (12.8,16)-- (14,16) node[right] {$A_{E}^{0,3}$};
     \draw[-latex] (16.8,16)-- (18,16) node[right] {$A_{E}^{0,4}$};
      \draw[-latex] (20.8,16)-- (22,16) node[right] {$0$};
      \draw[-latex] (3,17)-- (3,18.5);
\draw[-latex] (7,17)-- (7,18.5);
\draw[-latex] (11,17)-- (11,18.5);
\draw[-latex] (15,17)-- (15,18.5);
\draw[-latex] (19,17)-- (19,18.5);
\draw  (3.8,-3.3) node[left] {$0$}
           (7.8,-3.3) node[left] {$0$}
                      (11.8,-3.3) node[left] {$0$}
                                 (15.8,-3.3) node[left] {$0$}
                                            (19.8,-3.3) node[left] {$0$};
\draw  (1,0) node[left] {$0$}
         (1,4) node[left] {$0$}
            (1,8) node[left] {$0$}
              (1,12) node[left] {$0$}
               (1,16) node[left] {$0$};
   \draw  (3.8,19.3) node[left] {$0$}
           (7.8,19.3) node[left] {$0$}
                      (11.8,19.3) node[left] {$0$}
                                 (15.8,19.3) node[left] {$0$}
                                            (19.8,19.3) node[left] {$0$};
\draw[red] (1,7)--(1,17)-- (13,17)--(13,7)--(1,7);
 \draw[green] (5,3)--(5,13)-- (17,13)--(17,3)--(5,3);
    \draw[blue] (9,-1)--(9,11)-- (21,11)--(21,-1)--(9,-1);
           \end{tikzpicture}
\end{center}
We here state some facts which is necessary for the proof of Theorem \ref{Poisson-pjbd-formula}.
\begin{prop}\label{prop2}
Let $(X,\pi)$ be a holomorphic Poisson manifold.
Then for any $d$-closed forms $\alpha$ and
  $\beta$ on $X$, the bracket $[\alpha,\beta]_{\partial_{\pi}}$ is $d$-exact.
\end{prop}

\begin{proof}
Let $d_{\pi}:=\iota_{\pi}d-d\iota_{\pi}=\partial_{\pi}+\bar{\partial}_{\pi}$.
According to a result by Sharygin--Talalaev \cite[Lemma 5]{ST08},
the Lie bracket
$$
[\alpha,\beta]_{d_{\pi}}=
(-1)^k\big(d_{\pi}(\alpha\wedge \beta)-(d_{\pi}\alpha) \wedge \beta
               -(-1)^{k} \alpha\wedge (d_{\pi}\beta)\big),
               \alpha\in A_{X}^{k},\ \beta\in A_{X}^{l},
$$
associated to $d_{\pi}$ is $d$-exact
if both $\alpha$ and $\beta$ are $d$-closed.
Due to the Lemma \ref{lem1}, we obtain that $\bar{\partial}_{\pi}$ is zero,
and hence we get $d_{\pi}=\partial_{\pi}$.
Therefore we in fact have $[\alpha,\beta]_{d_{\pi}}=[\alpha,\beta]_{\partial_{\pi}}$,
and consequently the Proposition holds.
\end{proof}

We are now in a position to give the proof of Theorem \ref{Poisson-pjbd-formula}.
\begin{proof}[Proof of Theorem \ref{Poisson-pjbd-formula}]
Using the same notations as above, the morphism of the first pages
$\Psi_{1}:(\mathcal{G}^{\bullet,\bullet}_{1},d_{1}=\partial_{\pi})
\rightarrow(\mathcal{L}^{\bullet,\bullet}_{1},\bar{d}_{1}=\partial_{\tilde{\pi}})$
is explicitly expressed as:
\begin{equation}
\Psi_{1}=\sum_{i=0}^{c-1}h^{i}\wedge \rho^{\star}(-):
\bigoplus_{i=0}^{c-1} \bigl(H_{\bar{\partial}}^{\bullet, \bullet}(Z)[-c+1+i, -i], \partial_{\pi}\bigr)
\longrightarrow  \bigl(H_{\bar{\partial}}^{\bullet, \bullet}(E),\partial_{\tilde{\pi}}\bigr).
\end{equation}
We claim that $\Psi_{1}$ commutes with $\partial_{\pi}$ and $\partial_{\tilde{\pi}}$.
Note that if the $\partial\bar{\partial}$-lemma holds on $Z$, 
it also holds on $E$ (cf. \cite[Corollary 12]{ASTT20}).
Since $h$ is a $d$-closed real $(1,1)$-form on $E$, it follows from Proposition \ref{prop2} that
$[h,h]_{\partial_{\tilde{\pi}}}$ is $d$-exact.
On the one hand, since $\partial h=0$, we get
\begin{equation}\label{partial=0}
\partial([h,h]_{\partial_{\tilde{\pi}}})=
[\partial h,h]_{\partial_{\tilde{\pi}}}+[h,\partial h]_{\partial_{\tilde{\pi}}}=0.
\end{equation}
On the other hand, since $[h,h]_{\partial_{\tilde{\pi}}}$ is $d$-exact,
it follows from \eqref{partial=0} that
\begin{equation}\label{d=0}
0=d([h,h]_{\partial_{\tilde{\pi}}})=\bar{\partial}([h,h]_{\partial_{\tilde{\pi}}}).
\end{equation}
From the $\partial\bar{\partial}$-lemma, we obtain
$[h,h]_{\partial_{\tilde{\pi}}}=\partial\bar{\partial}\beta$ for some $\beta$ on $E$.
This implies that $[h,h]_{\partial_{\tilde{\pi}}}$ represents a zero class in $H^{1,2}_{\bar{\partial}}(E)$ and therefore we get
$$
0=[\partial_{\tilde{\pi}}h^{i}]\in H^{i-1,i}_{\bar{\partial}}(E).
$$
For any $[\alpha]\in H^{\bullet,\bullet}_{\bar{\partial}}(Z)$, we have
$\partial\alpha\in\ker\,\bar{\partial}\cap\mathrm{im\,\partial}$.
Since $Z$ satisfies the $\partial\bar{\partial}$-lemma,
 there exists $\xi$ on $Z$ such that $\partial\alpha=\partial\bar{\partial}\xi$.
Put $\tilde{\alpha}=\alpha-\bar{\partial}\xi$.
Then we get $[\alpha]=[\tilde{\alpha}]$ and $\partial\tilde{\alpha}=0$.
In what follows, we always choose the $\partial$-closed representatives
of the Dolbeault cohomology classes in $H^{\bullet,\bullet}_{\bar{\partial}}(Z)$.
Let $[\alpha]\in H^{\bullet,\bullet}_{\bar{\partial}}(Z)$, then we have
\begin{eqnarray*}
  \partial_{\tilde{\pi}}(\Psi_{1}([\alpha]))
  &=& \partial_{\tilde{\pi}}(\sum_{i=0}^{c-1}[h^{i}\wedge\rho^{\star}(\alpha)]))
   = \sum_{i=0}^{c-1}[\partial_{\tilde{\pi}}(h^{i}\wedge\rho^{\star}(\alpha))] \\
  &=& \sum_{i=0}^{c-1}[(\partial_{\tilde{\pi}}h^{i})\wedge\rho^{\star}(\alpha)+
       h^{i}\wedge(\partial_{\tilde{\pi}}\circ\rho^{\star}(\alpha))
       +[h^{i},\rho^{\star}(\alpha)]_{\partial_{\tilde{\pi}}}].
\end{eqnarray*}
Note that $\partial_{\tilde{\pi}}h^{i}$ is $\bar{\partial}$-exact and $\alpha$ is $\bar{\partial}$-closed.
We obtain that $[(\partial_{\tilde{\pi}}h^{i})\wedge\rho^{\star}(\alpha)]=0$.
Consider $\gamma:=[h^{i},\rho^{\star}(\alpha)]_{\partial_{\tilde{\pi}}}$.
From Proposition \ref{prop2}, we know that $\gamma$ is $d$-exact.
Notice that both $h$ and $\alpha$ are $\partial$-closed.
This implies
$\partial\gamma=[\partial h^{i},\rho^{\star}(\alpha)]_{\partial_{\tilde{\pi}}}+
[h^{i},\partial\rho^{\star}(\alpha)]_{\partial_{\tilde{\pi}}}=0$.
Furthermore, we get that $\gamma$ is $\bar{\partial}$-closed.
By
the $\partial\bar{\partial}$-lemma on $E$,
we get $\gamma=\partial\bar{\partial}\eta$ for some $\eta$ on $E$.
This implies $[\gamma]=0$ in the Dolbeault cohomology group.
Consequently, we are led to the conclusion
\begin{eqnarray*}
  \partial_{\tilde{\pi}}(\Psi_{1}([\alpha]))
  &=& \sum_{i=0}^{c-1}[h^{i}\wedge\partial_{\tilde{\pi}}(\rho^{\star}(\alpha))]
   = \sum_{i=0}^{c-1}[h^{i}\wedge\rho^{\star}(\partial_{\pi}(\alpha))]
   = \Psi_{1}(\partial_{\pi}([\alpha])).
\end{eqnarray*}
The morphism $\Psi_{1}$ induces the morphism between the second pages of the spectral sequences:
$$
\Psi_{2}=H(\Psi_{1}):\mathcal{G}^{\bullet,\bullet}_{2}\cong
H(\mathcal{G}^{\bullet,\bullet}_{1},d_{1})
\longrightarrow H(\mathcal{L}^{\bullet,\bullet}_{1},\bar{d}_{1})\cong
\mathcal{L}^{\bullet,\bullet}_{2}.
$$
Assume that
$\Psi_{r}:\mathcal{G}^{\bullet,\bullet}_{r}\longrightarrow\mathcal{L}^{\bullet,\bullet}_{r}$
commutes with the differentials $d_{r}$ and $\bar{d}_{r}$, where $r\geq2$.
Then we get the induced morphism of the $(r+1)$-pages:
$$
\Psi_{r+1}=H(\Psi_{r}):\mathcal{G}^{\bullet,\bullet}_{r+1}\cong
H(\mathcal{G}^{\bullet,\bullet}_{r},d_{r})
\longrightarrow H(\mathcal{L}^{\bullet,\bullet}_{r},\bar{d}_{r})\cong
\mathcal{L}^{\bullet,\bullet}_{r+1}.
$$
We claim that $\Psi_{r+1}$ commutes with $d_{r+1}$ and $\bar{d}_{r+1}$.
Suppose $\alpha\in A^{p,q}_{Z}$ represents a class in $\mathcal{G}^{p,q}_{r+1}$; equivalently, $\alpha$ is a cocycle in $\mathcal{G}^{p,q}_{j}$ for all $1\leq j\leq r$.
Without loss of generality, we assume that $\alpha$ is $\partial$-closed.
Since $\alpha$ lives to $\mathcal{G}^{p,q}_{r+1}$, by the $\partial\bar{\partial}$-lemma,
it can be extended to a zig-zag of length $r+1$ such that the tail is $\partial$-exact.
Denote the tail by $\partial\beta$ and then we get
$d_{r+1}([\alpha]_{r+1})=[\partial_{\pi}\partial\beta]_{r+1}$.
Observe that the form $\tilde{\alpha}:=h^{i}\wedge\rho^{\star}\alpha$
represents a class in $\mathcal{L}^{p,q}_{1}$.
Using the $\partial\bar{\partial}$-lemma again, the form $\tilde{\alpha}$
can be extended to a zig-zag of length $(r+1)$ which has the tail
$\tilde{\eta}:=h^{i}\wedge\rho^{\star}(\partial\beta)+\partial\tilde{\gamma}$,
where $\tilde{\gamma}$ is a form on $E$.
This implies that $\tilde{\alpha}$ lives to $\mathcal{L}^{p,q}_{r+1}$ and hence we get
$$
\bar{d}_{r+1}(\Psi_{r+1}([\alpha]_{r+1}))=
\bar{d}_{r+1}([\tilde{\alpha}]_{r+1})=[\partial_{\tilde{\pi}}\tilde{\alpha}]_{r+1}.
$$
Via a straightforward computation, we have
\begin{eqnarray*}
  \partial_{\tilde{\pi}}\tilde{\eta}
  &=&(\partial_{\tilde{\pi}}h^{i})\wedge\rho^{\star}(\partial\beta)+
    h^{i}\wedge\partial_{\tilde{\pi}}\rho^{\star}(\partial\beta)+
    [h^{i}, \rho^{\star}(\partial\beta)]_{\partial_{\tilde{\pi}}}+
    \partial_{\tilde{\pi}}\partial\tilde{\gamma}\\
  &=&  h^{i}\wedge\rho^{\star}(\partial_{\pi}\partial\beta)-
     \partial\bigl((\partial_{\tilde{\pi}}h^{i})\wedge\rho^{\star}\beta\bigr)-
     \partial([h^{i}, \rho^{\star}\beta]_{\partial_{\tilde{\pi}}})-
     \partial\partial_{\tilde{\pi}}\tilde{\gamma}\\
  &:=& h^{i}\wedge\rho^{\star}(\partial_{\pi}\partial\beta)+\partial\tilde{\zeta}.
\end{eqnarray*}
Note that both $\partial_{\tilde{\pi}}\tilde{\eta}$ and $\partial_{\pi}\partial\beta$
are $\bar{\partial}$-closed, and so is the form $\partial\tilde{\zeta}$.
Due to the $\partial\bar{\partial}$-lemma, we get
$\partial\tilde{\zeta}=\partial\bar{\partial}\tilde{\omega}$ for some $\tilde{\omega}$ on $E$.
This implies that $\partial_{\tilde{\pi}}\tilde{\eta}$ and $h^{i}\wedge\rho^{\star}(\partial_{\pi}\partial\beta)$
represent the same class in $\mathcal{L}_{1}$, i.e.,
$[\partial_{\tilde{\pi}}\tilde{\eta}]_{1}=
[h^{i}\wedge\rho^{\star}(\partial_{\pi}\partial\beta)]_{1}
$
and therefore we get
$$
\tilde{d}_{r+1}(\Psi_{r+1}([\alpha]_{r+1}))=
[\partial_{\tilde{\pi}}\tilde{\eta}]_{r+1}=
[h^{i}\wedge\rho^{\star}(\partial_{\pi}\partial\beta)]_{r+1}=
\Psi_{r+1}(d_{r+1}([\alpha]_{r+1})).
$$
Inductively, we obtain a morphism of spectral sequences
$$
\Psi_{r}:(\mathcal{G}^{\bullet,\bullet}_{r},d_{r})
\longrightarrow(\mathcal{L}^{\bullet,\bullet}_{r},\bar{d}_{r}).
$$
Thanks to the projective bundle formula for Dolbeault cohomology,
we get that $\Psi_{1}$ is an isomorphism.
By a result on the convergence of spectral sequences \cite[Theorem 3.4]{Mc01},
for all $r$,  $1< r\leq \infty$,
$\Psi_{r}:\mathcal{G}^{\bullet,\bullet}_{r}\rightarrow \mathcal{L}^{\bullet,\bullet}_{r}$
is an isomorphism, and this completes the proof.
\end{proof}

\begin{example}\label{p-x}
Suppose $(X,\pi_X)$ and $(\mathbb{P}^{n},\pi_{\mathbb{P}^{n}})$ are two holomorphic Poisson manifolds.
Then the product manifold $E=X\times \mathbb{P}^{n}$ can be thought of as the
projectivization of the trivial vector bundle $X\times\mathbb{C}^{n+1}\rightarrow X$.
In particular, if we view $E$ as the product of Poisson manifolds $X$ and $\mathbb{P}^{n}$,
then the associated product Poisson structure
$\tilde{\pi}_E$ on $E$ satisfying
that the two projections $\rho_1:E\rightarrow X$
and $\rho_2:E\rightarrow \mathbb{P}^{n}$ are Poisson maps with involution property
(cf. \cite[Proposition 2.5]{L-GPV13}).
Notice that
the first Chern class $[h_E]$ of the tautological line bundle over $E$
is the pullback of the first Chern class  $[h_{\mathbb{P}^{n}}]$ of the tautological line bundle over $\mathbb{P}^{n}$ via $\rho_2^\star$,
hence by the involution property of $\tilde{\pi}_E$, for any $\alpha\in\A^{\bullet,\bullet}_X$,
$$[h^i_E, \rho_1^\star(\alpha)]_{\partial_{\tilde{\pi}_E}}
=[\rho_2^\star(h^i_{\mathbb{P}^{n}}), \rho_1^\star(\alpha)]_{\partial_{\tilde{\pi}_E}}=0.$$
This means in this special case,
the map $\Psi$ is a well-defined morphism between the double complexes $G$ and $L$.
With the classical argument, $\Psi$ induces an isomorphism
$$
H_{k}(E,\tilde{\pi}_E)\cong H_{k-n}(X, \pi_X)^{\oplus (n+1)},
$$
for any integer $k\geq0$.
Especially, if $X$ is a point, then we have
$$
H_{k}(\mathbb{P}^{n},\pi_{\mathbb{P}^{n}})=
\left\{ \begin{array}{ll}
\mathbb{C}^{n+1},&k=n,\\
0, & k\neq n.
\end{array} \right.$$

Assume that $X=\mathbb{P}^{1}$ and $n=1$.
The moduli space of holomorphic Poisson structures on $E=\mathbb{P}^{1}\times \mathbb{P}^{1}$
is isomorphic to $\mathbb{C}^{9}$ (cf. \cite[Proposition 2.2]{HX11}).
Consider the projective bundle $\rho: E\rightarrow \mathbb{P}^{1}$.
Let $\pi_{E}$ be an arbitrary holomorphic Poisson structure on $E$.
The holomorphic Koszul--Brylinski homology of $(E,\pi_{E})$ has been computed by
Sti\'{e}non \cite[Theorem 7.2]{Sti11}, and here we present a new proof of his result by
applying the projective bundle formula (Theorem \ref{Poisson-pjbd-formula}).
For the dimension reason, each holomorphic bivector field on $\mathbb{P}^{1}$ is zero.
Thus we get $\rho_{\star}(\pi_E)=0$.
As a corollary, we obtain
$$
H_{k}(E,\pi_{E})\cong H_{k-1}(\mathbb{P}^{1}, \pi=0)^{\oplus 2}
\cong
\Big[\bigoplus_{p-q=2-k} H^{q}(\mathbb{P}^{1}, \Omega_{\mathbb{P}^{1}}^{p})\Big]^{\oplus 2}.
$$
A straightforward computation shows $H_{k}(E, \pi_{E})=0$ when $k=0, 1,3,4$.
When $k=2$, we have
$$
H_{2}(E,\pi_{E})\cong H_{1}(\mathbb{P}^{1}, \pi=0)^{\oplus 2}\cong \mathbb{C}^{4},
$$
since $H_{1}(\mathbb{P}^{1}, \pi=0)\cong
H_{\bar{\partial}}^{0,0}(\mathbb{P}^{1})\oplus H_{\bar{\partial}}^{1,1}(\mathbb{P}^{1})\cong \mathbb{C}^{2}$.
\end{example}


\section{Comparison under Poisson blow-ups}\label{Comparison under Poisson blow-ups}

The main purpose of this section is to prove the blow-up formula for holomorphic Koszul--Brylinski homology of compact holomorphic Poisson manifolds.

\subsection{Relative Koszul--Brylinski homology}
Let $X$ be a compact complex manifold and $\jmath:Z\hookrightarrow X$ a closed complex submanifold of codimension $c$.
Consider the natural morphism
$$
\jmath^{\star}:\Omega_{X}^{s}\longrightarrow \jmath_{\ast}\Omega^{s}_{Z}
$$
which is defined as follows:
\begin{eqnarray*}
\jmath^{\star}(V):\Gamma(V,\Omega_{X}^{s})
&\longrightarrow& \Gamma(V,\jmath_{\ast}\Omega^{s}_{Z}) \\
\alpha &\longmapsto& (\jmath_{V\cap Z})^{\star} \alpha,
\end{eqnarray*}
where $V\subset X$ is an open subset and $\jmath_{V\cap Z}:V\cap Z\hookrightarrow V$
is the holomorphic inclusion.
We also define the sheaf morphism
$$
\jmath^{\star}: \A_{X}^{p,q} \longrightarrow \jmath_{\ast}\A_{Z}^{p,q}.
$$
in a similar way.

\begin{defn}[\cite{RYY20,YY20}]
The kernel sheaves
$$
\K_{X,Z}^{s}:=\ker\big(\jmath^{\star}:\Omega_{X}^{s}\rightarrow \jmath_{\ast}\Omega^{s}_{Z}\big)
$$
and
$$
\K_{X,Z}^{p,q}:=\ker\big(\jmath^{\star}:\A_{X}^{p,q} \rightarrow \jmath_{\ast}\A_{Z}^{p,q}\big)
$$
are called the {\it $s$-th relative Dolbeault sheaf} and {\it $(p,q)$-th
relative Dolbeault sheaf} with respect to $Z$.
\end{defn}

There exist two natural short exact sequences
\begin{equation}\label{rel-pureDolsh}
\xymatrix@C=0.5cm{
0 \ar[r]^{} & \K_{X,Z}^{s} \ar[r]^{} & \Omega_{X}^{s} \ar[r]^{\jmath^{\star}}
& \jmath_{\ast}\Omega_{Z}^{s}\ar[r]^{} & 0}
\end{equation}
and
\begin{equation}\label{rel-Dolsh}
\xymatrix@C=0.5cm{
0 \ar[r]^{} & \K_{X,Z}^{p, q} \ar[r]^{} & \A_{X}^{p, q} \ar[r]^{\jmath^{\star}}
& \jmath_{\ast}\A_{Z}^{p, q}\ar[r]^{} & 0,}
\end{equation}
where $\K_{X,Z}^{p,\bullet}$ is a fine resolution of $\K_{X,Z}^{p}$.
Consider the holomorphic Poisson manifold $(X,\pi)$ together
with the holomorphic Poisson submanifold $(Z,\pi|_{Z})$.
Then we have

\begin{lem}
There exists a short exact sequence of sheaf complexes on $X$:
\begin{equation}\label{XZ-sh-exact-sequecne}
\xymatrix@C=0.5cm{
  0 \ar[r] & (\K^{\bullet}_{X,Z}, \partial_{\pi}) \ar[r]^{} & (\Omega_{X}^{\bullet},
  \partial_{\pi}) \ar[r]^{\jmath^{\star}\;\;\;\;\;} & \jmath_{\ast}(\Omega_{Z}^{\bullet},
  \partial_{\pi|_{Z}})[-c] \ar[r] & 0,}
\end{equation}
where $c=\mathrm{codim}_{\mathbb{C}}Z$.
\end{lem}

\begin{proof}
Since $\jmath^{\star} \circ \iota_{\pi}=\iota_{\pi|_{Z}}\circ \jmath^{\star}$
and
$\jmath^{\star} \circ \partial=\partial\circ \jmath^{\star}$,
we have
$\jmath^{\star} \circ \partial_{\pi}=\partial_{\pi|_{Z}}\circ \jmath^{\star}$.
Hence, for any $k\in \mathbb{Z}$,
there is a well-defined induced operator $\partial_{\pi}: \K^{k}_{X,Z}\rightarrow \K^{k-1}_{X,Z}$,
i.e., $(\K^{\bullet}_{X,Z}, \partial_{\pi})$ is a well-defined sheaf complex.
Moreover, by the short exact sequence \eqref{rel-pureDolsh},
there is a commutative diagram of short exact sequences of sheaves
\begin{equation*}
\xymatrix@C=0.5cm{
 0 \ar[r]^{} &  \K_{X,Z}^{k} \ar[d]_{\partial_{\pi}} \ar[r]^{} & \Omega_{X}^{k} \ar[d]_{\partial_{\pi}} \ar[r]^{\jmath^{\star}} & \jmath_{\ast} \Omega_{Z}^{k} \ar[d]_{\partial_{\pi|_{Z}}} \ar[r]^{} & 0 \\
0 \ar[r] &  \K_{X,Z}^{k-1} \ar[r]^{} &  \Omega_{X}^{k-1} \ar[r]^{\jmath^{\star}}
&  \jmath_{\ast}\Omega_{Z}^{k-1} \ar[r] &0. }
\end{equation*}
The lemma follows immediately from the above commutative diagram.
\end{proof}

\begin{defn}
The sheaf complex $(\K_{X,Z}^{\bullet}, \partial_{\pi})$ is called the
{\it relative Koszul--Brylinski complex} of $(X,\pi)$ with respect to $(Z,\pi|_{Z})$,
and its $k$-th hypercohomology
$$
H_{k}(X,Z;\pi):=\mathbb{H}^{k}(X, (\K_{X,Z}^{\bullet}, \partial_{\pi}))
$$
is called the $k$-th {\it relative Koszul--Brylinski homology} of $(X,\pi)$ with respect to  $(Z,\pi|_{Z})$.
\end{defn}

Taking the hypercohomology of the short exact sequence \eqref{XZ-sh-exact-sequecne},
we get a long exact sequence:
\begin{equation*}
\xymatrix@C=0.4cm{
\cdots \ar[r]^{}
& H_{k}(X, Z;\pi)\ar[r]^{}
& H_{k}(X, \pi) \ar[r]^{}
& H_{k-c}(Z, \pi|_{Z})  \ar[r]^{}
& H_{k+1}(X,Z;\pi)  \ar[r]^{} & \cdots}
\end{equation*}
Similarly to Proposition \ref{KB-equal1}, we have the following:

\begin{prop}
The relative Koszul--Brylinski complex $(\K_{X,Z}^{\bullet}, \partial_{\pi})$
is quasi-isomorphic to
the total sheaf complex of the double sheaf complex $(\K_{X,Z}^{\bullet,\bullet},\partial_{\pi},\bar{\partial})$.
\end{prop}

\begin{proof}
As we mentioned above, for each $p\in \mathbb{Z}$, the relative
Dolbeault sheaf $\K_{X,Z}^{p}$ admits a fine resolution $\K_{X,Z}^{p,\bullet}$.
Then the proposition follows.
\end{proof}


\subsection{Proof of Theorem \ref{main-thm}}
Given a compact  holomorphic Poisson manifold $(X, \pi)$ with a holomorphic Poisson closed submanifold
$\jmath:(Z,\pi|_{Z}) \hookrightarrow (X,\pi)$
with complex codimension $c\geq 2$.
Let $\varphi: \tilde{X}\rightarrow X$ be the blow-up of $X$ along $Z$
with exceptional divisor $E:=\varphi^{-1}(Z)$.
Suppose that the transverse Poisson structure $\pi_{N}=0$.
Due to Proposition \ref{Constr-Poisson-bp},
we get a commutative diagram for the blow-up in the holomorphic Poisson category:
\begin{equation}\label{poiss-bl-dia}
\vcenter{
\xymatrixcolsep{3pc}
\xymatrix{
(E,\tilde{\pi}|_{E}) \ar[d]_{\rho=\varphi|_E}
\ar@{^{(}->}[r]^{\tilde{\jmath}} & (\tilde{X},\tilde{\pi})\ar[d]^{\varphi}\\
(Z,\pi|_{Z}) \ar@{^{(}->}[r]^{\jmath} & (X,\pi).
}}
\end{equation}
The following lemma plays a crucial role in the proof Theorem \ref{main-thm}.

\begin{lem}\label{Poisson-kersheaf-iso}
The pullback $\varphi^{\star}$ naturally induces an isomorphism
$$
\varphi^{\star}:
H_{k}(X,Z;\pi)
\stackrel{\simeq}\longrightarrow
H_{k}(\tilde{X},E;\tilde{\pi}).
$$
for any $k\in \mathbb{Z}$.
\end{lem}

\begin{proof}
Note that there exists a natural morphism of bounded double complexes
$$
\varphi^{\star}:
\bigl(\Gamma(X,\K_{X,Z}^{\bullet,\bullet}),\partial_{\pi}, \bar{\partial})\bigr)
\longrightarrow
\bigl(\Gamma(\tilde{X},\K_{\tilde{X},E}^{\bullet,\bullet}),\partial_{\tilde{\pi}}, \bar{\partial})\bigr).
$$
Furthermore, we have two spectral sequences:
\begin{itemize}
\item $\{E_{r},d_{r}\}$,
associated to $(\Gamma(X,\K_{X,Z}^{\bullet,\bullet}),\partial_{\pi}, \bar{\partial}))$ ,
converges to the relative Koszul--Brylinski homology $H_{\bullet}(X,Z;\pi)$
with the $E_{1}$-page given by
\begin{equation*}
E_{1}^{p,q}= H^{q}(X, \K^{n-p}_{X,Z});
\end{equation*}
\item $\{\tilde{E}_{r},\tilde{d}_{r}\}$,
associated to $(\Gamma(\tilde{X},\K_{\tilde{X},E}^{\bullet,\bullet}),
\partial_{\tilde{\pi}}, \bar{\partial}))$, converges to the relative Koszul--Brylinski homology
$H_{k}(\tilde{X},E;\tilde{\pi})$ with the $\tilde{E}_{1}$-page given by
\begin{equation*}
\tilde{E}_{1}^{p,q}=H^{q}(\tilde{X},\K_{\tilde{X},E}^{n-p}).
\end{equation*}
\end{itemize}
For any $r\geq1$, the morphism $\varphi^{\star}$
induces a morphism of the spectral sequences
$$
\varphi^{\star}_{r}:E_{r}\longrightarrow\tilde{E}_{r}
$$
and hence a morphism of relative Koszul--Brylinski homologies
$$
\varphi^{\star}:
H_{k}(X,Z;\pi)
\longrightarrow
H_{k}(\tilde{X},E;\tilde{\pi}).
$$
Due to \cite[Lemma 4.5]{RYY20},
the pullback of differential forms induces an isomorphism
\begin{equation*}
\varphi^{\star}:
H^{q}(X, \K_{X,Z}^{n-p})
\longrightarrow
H^{q}(\tilde{X}, \K_{\tilde{X},E}^{n-p}),
\end{equation*}
for any $0\leq p, q\leq n$.
It follows that $\varphi^{\star}_{1}:E_{1}\rightarrow\tilde{E}_{1}$ is an isomorphism.
Consequently, by the standard result in the spectral sequence theory, we get that
$\varphi^{\star}_{r}$ is isomorphic for any $r>1$ and therefore the assertion holds.
\end{proof}

Now we are in the position to prove the blow-up formula of holomorphic Koszul--Brylinski homology.

\begin{proof}[Proof of Theorem \ref{main-thm}]

For the pairs of compact holomorphic Poisson manifolds
$(X, Z)$ and $(\tilde{X},E)$,
we have two short exact sequences of sheaf complexes
\begin{equation}\label{Z-holomor-seq}
\xymatrix@C=0.5cm{
  0 \ar[r] & (\K_{X,Z}^{\bullet},\partial_{\pi}) \ar[r]^{}
  & (\Omega_{X}^{\bullet}, \partial_{\pi})  \ar[r]^{ \jmath^{\star}\;\;\;\;}
  &  \jmath_{\ast}(\Omega_{Z}^{\bullet},\partial_{\pi|_{Z}})[-c]  \ar[r] & 0}
\end{equation}
and
\begin{equation}\label{E-holomor-seq}
\xymatrix@C=0.5cm{
  0 \ar[r] & (\K_{\tilde{X},E}^{\bullet},\partial_{\tilde{\pi}}) \ar[r]^{}
  & (\Omega_{\tilde{X}}^{\bullet}, \partial_{\tilde{\pi}})
  \ar[r]^{\tilde{\jmath}^{\star}\;\;\;\;} & \tilde{\jmath}_{\ast}(\Omega_{E}^{\bullet},
  \partial_{\tilde{\pi}|_{E}})[-1]  \ar[r] & 0}.
\end{equation}
We next establish a commutative diagram of long exact sequences of
holomorphic Koszul--Brylinski homology
associated to \eqref{Z-holomor-seq} and \eqref{E-holomor-seq}.
Observe that each complex in \eqref{Z-holomor-seq} and
\eqref{E-holomor-seq} admits a natural fine resolution.
For the pair $(X, Z)$,
by the short exact sequence \eqref{rel-Dolsh},
we have the two commutative diagrams of short exact sequences
of fine sheaves:
\begin{equation*}
\xymatrix@C=0.5cm{
 0 \ar[r]^{} &  \K_{X,Z}^{p,q} \ar[d]_{\bar{\partial}} \ar[r]^{} & \A_{X}^{p,q}
 \ar[d]_{\bar{\partial}} \ar[r]^{\jmath^{\star}} & \jmath_{\ast} \A_{Z}^{p,q}
 \ar[d]_{\bar{\partial}} \ar[r]^{} & 0 \\
0 \ar[r] &  \K_{X,Z}^{p,q+1} \ar[r]^{} &  \A_{X}^{p,q+1} \ar[r]^{\jmath^{\star}}
&  \jmath_{\ast}\A_{Z}^{p,q+1} \ar[r] &0}
\end{equation*}
and
\begin{equation*}
\xymatrix@C=0.5cm{
 0 \ar[r]^{} &  \K_{X,Z}^{p,q} \ar[d]_{\partial_{\pi}} \ar[r]^{} & \A_{X}^{p,q}
 \ar[d]_{\partial_{\pi}} \ar[r]^{\jmath^{\star}} & \jmath_{\ast} \A_{Z}^{p,q}
 \ar[d]_{\partial_{\pi|_{Z}}} \ar[r]^{} & 0 \\
0 \ar[r] &  \K_{X,Z}^{p-1,q} \ar[r]^{} &  \A_{X}^{p-1,q} \ar[r]^{\jmath^{\star}}
&  \jmath_{\ast}\A_{Z}^{p-1,q} \ar[r] &0. }
\end{equation*}
As a consequence, there is a short exact sequence of \emph{double complexes}:
\begin{equation*}
\xymatrix@C=0.5cm{
 0
 \ar[r]^{} &
 \big(\Gamma(X, \K_{X,Z}^{\bullet,\bullet}), \partial_{\pi}, \bar{\partial}\big)
 \ar[r]^{} &
 \big(\Gamma(X, \A_{X}^{\bullet,\bullet}), \partial_{\pi}, \bar{\partial}\big)
 \ar[r]^{\jmath^{\star}\;\;\;\;} &
\big(\Gamma(Z, \A_{Z}^{\bullet,\bullet})[-c,0], \partial_{\pi|_{E}}, \bar{\partial}\big)
 \ar[r]^{}
 & 0.}
\end{equation*}
Similarly, for the pair $(\tilde{X},E)$ we have a short exact sequence:
\begin{equation*}
\xymatrix@C=0.5cm{
 0
 \ar[r]^{} &
 \big(\Gamma(\tilde{X}, \K_{\tilde{X},E}^{\bullet,\bullet}), \partial_{\tilde{\pi}}, \bar{\partial}\big)
 \ar[r]^{} &
 \big(\Gamma(\tilde{X}, \A_{\tilde{X}}^{\bullet,\bullet}), \partial_{\tilde{\pi}}, \bar{\partial}\big)
 \ar[r]^{\tilde{\jmath}^{\star}\;\;\;\;} &
\big(\Gamma(E, \A_{E}^{\bullet,\bullet})[-1,0], \partial_{\tilde{\pi}|_{E}}, \bar{\partial}\big)
 \ar[r]^{}
 & 0.}
\end{equation*}
Moreover,
the blow-up morphism $\varphi$ naturally induces a commutative diagram of short exact sequences:
\begin{equation*}
 \begin{tikzpicture}[scale=0.5]
\draw  (1.3,0) node[left] {$0$};
            \draw[-latex] (1.3,0)-- (3.5,0) node[right] {$\Gamma(X,\K_{X,Z}^{p,q})$};
  \draw[-latex] (8,0)-- (12,0) node[right] {$\Gamma(X,\A_{X}^{p,q})$};
   \draw[-latex] (16,0)-- (20,0) node[right] {$\Gamma(Z,\A_{Z}^{p,q})$};
   \draw[-latex] (24,0)-- (26.2,0) node[right] {$0$};
\draw  (5.3,-3) node[left] {$0$};
            \draw[-latex] (5.3,-3)-- (7.5,-3) node[right] {$\Gamma(X,\K_{X,Z}^{p-1,q})$};
  \draw[-latex] (12.3,-3)-- (16,-3) node[right] {$\Gamma(X,\A_{X}^{p-1,q})$};
   \draw[-latex] (20.7,-3)-- (24,-3) node[right] {$\Gamma(Z,\A_{Z}^{p-1,q})$};
   \draw[-latex] (28.7,-3)-- (31.2,-3) node[right] {$0$};
\draw  (1.3,-6) node[left] {$0$};
            \draw[-latex] (1.3,-6)-- (3.5,-6) node[right] {$\Gamma(\tilde{X},\K_{\tilde{X},E}^{p,q})$};
  \draw[-latex,dashed] (8,-6)-- (12,-6) node[right] {$\Gamma(\tilde{X},\A_{\tilde{X}}^{p,q})$};
   \draw[-latex,dashed] (16,-6)-- (20,-6) node[right] {$\Gamma(E,\A_{E}^{p,q})$};
   \draw[-latex] (24,-6)-- (26.2,-6) node[right] {$0$};
\draw  (5.3,-9) node[left] {$0$};
            \draw[-latex] (5.3,-9)-- (7.5,-9) node[right] {$\Gamma(\tilde{X},\K_{\tilde{X},E}^{p-1,q})$};
  \draw[-latex] (12.3,-9)-- (16,-9) node[right] {$\Gamma(\tilde{X},\A_{\tilde{X}}^{p-1,q})$};
   \draw[-latex] (20.7,-9)-- (24,-9) node[right] {$\Gamma(E,\A_{E}^{p-1,q})$};
   \draw[-latex] (28.7,-9)-- (31.2,-9) node[right] {$0$};
     \draw  (18.7,0.5) node[left] {$\jmath^{\star}$};
          \draw  (18.7,-5.5) node[left] {$\tilde{\jmath}^{\star}$};
               \draw  (22.7,-2.5) node[left] {$\jmath^{\star}$};
                        \draw  (22.7,-8.5) node[left] {$\tilde{\jmath}^{\star}$};
 \draw[-latex,red] (6,-0.5)-- (9,-2.5);
  \draw[-latex,red] (14,-0.5)-- (17,-2.5);
   \draw[-latex,red] (22,-0.5)-- (25,-2.5);
    \draw[-latex,red] (6,-6.5)-- (9.4,-8.5);
  \draw[-latex,red] (14,-6.5)-- (17,-8.5);
   \draw[-latex,red] (22,-6.5)-- (25,-8.5);
   \draw  (8.8,-1.3) node[left,red] {$\partial_\pi$};
      \draw  (16.8,-1.3) node[left,red] {$\partial_\pi$};
         \draw  (24.8,-1.3) node[left,red] {$\partial_\pi$};
            \draw  (9,-7.3) node[left,red] {$\partial_{\tilde{\pi}}$};
               \draw  (16.8,-7.3) node[left,red] {$\partial_{\tilde{\pi}}$};
                  \draw  (24.8,-7.3) node[left,red] {$\partial_{\tilde{\pi}}$};
 \draw[-latex,blue] (5.5,-0.5)-- (5.5,-5);
  \draw[-latex,blue,dashed] (13.5,-0.5)-- (13.5,-5);
   \draw[-latex,blue,dashed] (21.5,-0.5)-- (21.5,-5);
    \draw[-latex,blue] (9.5,-3.5)-- (9.5,-8);
  \draw[-latex,blue] (17.5,-3.5)-- (17.5,-8);
   \draw[-latex,blue] (25.5,-3.5)-- (25.5,-8);
      \draw  (9.7,-4.5) node[left,blue] {$\varphi^{\star}$};
        \draw  (17.7,-4.5) node[left,blue] {$\varphi^{\star}$};
          \draw  (25.7,-4.5) node[left,blue] {$\rho^{\star}$};
             \draw  (5.7,-1.5) node[left,blue] {$\varphi^{\star}$};
                \draw  (13.7,-1.5) node[left,blue] {$\varphi^{\star}$};
                   \draw  (21.7,-1.5) node[left,blue] {$\rho^{\star}$};
          \end{tikzpicture}
\end{equation*}
Therefore, we obtain a commutative diagram of double complexes
\begin{equation}\label{cube-diagram}
\vcenter{
\xymatrix@C=0.3cm{
 0
 \ar[r]^{} &
 \big(\Gamma(X, \K_{X,Z}^{\bullet,\bullet}), \partial_{\pi}, \bar{\partial}\big)
 \ar[d]_{\varphi^{\star}} \ar[r]^{} &
 \big(\Gamma(X, \A_{X}^{\bullet,\bullet}), \partial_{\pi}, \bar{\partial}\big)
  \ar[d]_{\varphi^{\star}} \ar[r]^{\jmath^{\star}\;\;\;\;} &
\big(\Gamma(Z, \A_{Z}^{\bullet,\bullet})[-c,0], \partial_{\pi|_{E}}, \bar{\partial}\big)
  \ar[d]_{\rho^{\star}}  \ar[r]^{} & 0 \\
 0
 \ar[r]^{} &
 \big(\Gamma(\tilde{X}, \K_{\tilde{X},E}^{\bullet,\bullet}), \partial_{\tilde{\pi}}, \bar{\partial}\big)
 \ar[r]^{} &
 \big(\Gamma(\tilde{X}, \A_{\tilde{X}}^{\bullet,\bullet}), \partial_{\tilde{\pi}}, \bar{\partial}\big)
 \ar[r]^{\jmath^{\star}\;\;\;\;} &
\big(\Gamma(E, \A_{E}^{\bullet,\bullet})[-1,0], \partial_{\tilde{\pi}|_{E}}, \bar{\partial}\big)
 \ar[r]^{}
 & 0.}}
\end{equation}
The commutative diagram \eqref{cube-diagram} above yields a commutative diagram
of long exact sequences of Koszul--Brylinski homologies:
\begin{equation}\label{KB-homology-long-seq}
\vcenter{
\xymatrix@C=0.5cm{
  \cdots \ar[r]^{}
  & H_{k}(X,Z;\pi)\ar[d]_{\varphi^{\ast}} \ar[r]^{}
  & H_{k}(X, \pi)\ar[d]_{\varphi^{\star}} \ar[r]^{}
  & H_{k-c}(Z, \pi|_{Z}) \ar[d]_{\rho^{\star}} \ar[r]^{}
  & H_{k+1}(X,Z;\pi) \ar[d]_{\varphi^{\star}} \ar[r]^{} & \cdots \\
   \cdots \ar[r]
  & H_{k}(\tilde{X},E;\tilde{\pi})\ar[r]^{}
  & H_{k}(\tilde{X}, \tilde{\pi}) \ar[r]^{}
  & H_{k-1}(E, \tilde{\pi}|_{E}) \ar[r]
  & H_{k+1}(\tilde{X},E;\tilde{\pi})\ar[r]&\cdots.}}
\end{equation}
In the above diagram,
by Lemma \ref{Poisson-kersheaf-iso},
for any $l\in \mathbb{Z}$,
the morphism
$$
\varphi^{\star}:
H_{l}(X,Z;\pi)
\stackrel{\simeq}\longrightarrow
H_{l}(\tilde{X},E;\tilde{\pi})
$$
is an isomorphism.
Moreover,
by Theorem \ref{KB-injective},
both the second and third vertical arrows in \eqref{KB-homology-long-seq} are injective.
Finally, by a standard diagram-chasing in \eqref{KB-homology-long-seq},
we get the following isomorphisms of finite dimensional $\mathbb{C}$-vector spaces:
$$
H_{k}(\tilde{X}, \tilde{\pi})
\cong H_{k}(X, \pi)\oplus \Big( H_{k-1}(E, \tilde{\pi}|_{E})/ \rho^{\star}H_{k-c}(Z, \pi|_{Z})\Big)
$$
Furthermore,
if $Z$ satisfies the $\partial\bar{\partial}$-lemma,
then by Theorem \ref{Poisson-pjbd-formula},
\begin{eqnarray*}
H_{k}(\tilde{X}, \tilde{\pi})
\cong H_{k}(X, \pi)\oplus  H_{k-c}(Z, \pi|_{Z})^{\oplus c-1}.
\end{eqnarray*}
This completes the proof of Theorem \ref{main-thm}.
\end{proof}


\subsection{Degeneracy of the Dolbeault--Koszul--Brylinski spectral sequence}

Let $(X, \pi)$ be a holomorphic Poisson manifold of complex dimension $n$.
Consider the Koszul--Brylinski double complex
$(\Gamma(X, \mathcal{A}_{X}^{\bullet,\bullet}), \partial_{\pi},\bar{\partial})$.
Inspired by the Fr\"{o}hlicher (or Hodge--de Rham) spectral sequence of complex manifolds,
we introduce the following:

\begin{defn}\label{Dol-Poi-spectral-seq}
The Fr\"{o}hlicher-type spectral sequence associated to the double complex
$(\Gamma(X, \mathcal{A}_{X}^{\bullet,\bullet}), \partial_{\pi},\bar{\partial})$ satisfying
\begin{equation}\label{Dol-Poisson-spec.-seq}
E_{1}^{s,t}:=H_{\bar{\partial}}^{n-s,t}(X)
\Longrightarrow
H_{n-s+t}(X,\pi),
\end{equation}
is called the {\it Dolbeault--Koszul--Brylinski spectral sequence of $(X, \pi)$}.
\end{defn}

As mentioned before,
for a compact complex manifold $X$ with the trivial holomorphic Poisson structure $\pi$,
the Dolbeault--Koszul--Brylinski spectral sequence degenerates at $E_{1}$-page and we have
$$
H_{k}(X, \pi=0)
\cong
\bigoplus\limits_{p-q=n-k}H^{q}(X, \Omega_{X}^{p})
\cong
\bigoplus\limits_{p-q=n-k}H^{p,q}_{\bar{\partial}}(X).
$$
Analogously to the Hodge--de Rham spectral sequence,
in general, the Dolbeault--Koszul--Brylinski spectral sequence \eqref{Dol-Poisson-spec.-seq}
does not degenerate at $E_{1}$-page (see for example in \S\S\ref{ex1-non-degeeracy-E1}).
We have the following result.

\begin{lem}
Let $(X, \pi)$ be a compact holomorphic Poisson manifold of complex dimension $n$.
Then its Dolbeault--Koszul--Brylinski spectral sequence degenerates at $E_{1}$-page if and only if
$$
\sum_{p-q=n-k} \dim_{\mathbb{C}}\, H_{\bar{\partial}}^{p,q}(X)=\dim_{\mathbb{C}} H_{k}(X,\pi),
$$
for any $0\leq k\leq 2n$.
\end{lem}

\begin{proof}
Observe that for a compact holomorphic Poisson manifold the holomorphic
Koszul--Brylinski homology groups are finite-dimensional;
moreover, the following inequality holds
$$
\dim_{\mathbb{C}}\, H_{k}(X,\pi)
\leq
\sum_{p-q=n-k} \dim_{\mathbb{C}}\, H_{\bar{\partial}}^{p,q}(X)
$$
for any $0\leq k\leq 2n$.
By definition, the $E_{1}$-degeneracy of the Dolbeault--Koszul--Brylinski spectral sequence is
equivalent to the condition
$$
\dim_{\mathbb{C}}\, H_{k}(X,\pi)
=
\sum_{p-q=n-k} \dim_{\mathbb{C}}\, H_{\bar{\partial}}^{p,q}(X),
$$
for any $0\leq k\leq 2n$.
\end{proof}

We are ready to give the proof of Theorem \ref{main-thm-2}.

\begin{proof}[Proof of Theorem \ref{main-thm-2}]
By the blow-up formula for Dolbeault cohomology \cite[Theorem 1.2]{RYY19},
we have
$$
\sum_{p-q=n-k} \dim_{\mathbb{C}}\, H_{\bar{\partial}}^{p,q}(\tilde{X})
=
\sum_{p-q=n-k} \biggl[ \dim_{\mathbb{C}}\, H_{\bar{\partial}}^{p,q}(X)
+\sum_{i=1}^{c-1} \dim_{\mathbb{C}}\, H_{\bar{\partial}}^{p-i,q-i}(Z)\biggr].
$$
Consequently, by Theorem \ref{main-thm}, we get
\begin{eqnarray*}
&&\dim_{\mathbb{C}}\, H_{k}(\tilde{X}, \tilde{\pi})-\sum_{p-q=n-k}
\dim_{\mathbb{C}}\, H_{\bar{\partial}}^{p,q}(\tilde{X}) \\
&=&
\biggl[\dim_{\mathbb{C}}\, H_{k}(X,\pi)-\sum_{p-q=n-k}
\dim_{\mathbb{C}}\, H_{\bar{\partial}}^{p,q}(X)\biggr] \\
&&
+(c-1)\dim_{\mathbb{C}}\, H_{k-c}(Z,\pi|_{Z})-
\sum_{p-q=n-k}\biggl[\sum_{i=1}^{c-1} \dim_{\mathbb{C}}\, H_{\bar{\partial}}^{p-i,q-i}(Z)\biggr]\\
&=&
\biggl[\dim_{\mathbb{C}}\, H_{k}(X,\pi)-\sum_{p-q=n-k}
\dim_{\mathbb{C}}\, H_{\bar{\partial}}^{p,q}(X)\biggr] \\
& &+
(c-1)\biggl[\dim_{\mathbb{C}}\, H_{k-c}(Z,\pi|_{Z})-
\sum_{s-t=(n-c)-(k-c)}\dim_{\mathbb{C}}\, H_{\bar{\partial}}^{s,t}(Z)\biggr]\
\end{eqnarray*}
for $0\leq k\leq 2n$.
If the Dolbeault--Koszul--Brylinski spectral sequence degenerates at $E_{1}$-pages
for $(X, \pi)$ and $(Z, \pi_{Z})$,
then it immediately follows that the Dolbeault--Koszul--Brylinski spectral sequence
degenerates at $E_{1}$-pages for $(\tilde{X}, \tilde{\pi})$.
Conversely,
if the Dolbeault--Koszul--Brylinski spectral sequence
degenerates at $E_{1}$-pages for $(\tilde{X}, \tilde{\pi})$, then we obtain the following equalities
\begin{eqnarray*}
0&=&\underbrace{\dim_{\mathbb{C}}\, H_{k}(X,\pi)-\sum_{p-q=n-k}
\dim_{\mathbb{C}}\, H_{\bar{\partial}}^{p,q}(X)}_{\leq 0} \\
&+&
\underbrace{(c-1)\biggl[\dim_{\mathbb{C}}\, H_{k-c}(Z,\pi|_{Z})-
\sum_{s-t=(n-c)-(k-c)}\dim_{\mathbb{C}}\, H_{\bar{\partial}}^{s,t}(Z)\biggr],}_{\leq 0}
\end{eqnarray*}
which implies that the $E_{1}$-degeneracy holds for $(X, \pi)$ and $(Z,\pi|_{Z})$.
\end{proof}


\section{Examples}\label{examples}
In this section, as applications of the main theorems, we compute the
Koszul--Brylinski homology for some special holomorphic Poisson manifolds,
such as del Pezzo surfaces and Iwasawa manifolds.

\subsection{del Pezzo surfaces}

Recall that a {\it del Pezzo surface} is a smooth Fano surface which
is exactly one of the following: $\mathbb{P}^{1}\times \mathbb{P}^{1}$,
$\mathbb{P}^{2}$ and blow-up of $\mathbb{P}^{2}$ at $r$ $(1\leq r \leq 8)$
generic points (denoted by $M_r$).
The holomorphic Koszul--Brylinski homology of $\mathbb{P}^{1}\times \mathbb{P}^{1}$
has been computed in Example \ref{p-x}; see also {\cite[Theorem 7.2]{Sti11}}.
We now consider the rest cases.
Define the space
$$
V_{r}^{2}=\{ \textup{holomorphic bi-vector fields on $\mathbb{P}^{2}$
vanishing at the blow-up points of $M_{r}$} \}.
$$
By a result of Kodaira \cite[page 225]{Kod86},
the blow-up transformation $\varphi: M_{r}\rightarrow\mathbb{P}^{2}$
induces an isomorphism from the space of holomorphic bi-vector fields on $M_{r}$ to the space $V_{r}^{2}$.
Equivalently, the holomorphic Poisson structures $\pi$ on $\mathbb{P}^{2}$
vanishing at the blow-up points of $M_{r}$ are one-one corresponding to
the holomorphic Poisson structures $\tilde{\pi}$ on $M_{r}$ such that $\varphi$ is a Poisson morphism.

In general, given a holomorphic Poisson structure $\pi$ on $\mathbb{P}^{n}$, the $E_1$-page of the
Dolbeault--Koszul--Brylinski spectral sequence of $(\mathbb{P}^{n}, \pi)$ is
$$
E_{1}^{s,t}=H^{t}(\mathbb{P}^{n},\Omega^{n-s})=
  \left\{ \begin{array}{ll}
\mathbb{C},&s+t=n,\\
0, & \textrm{otherwise}.
\end{array} \right.
$$
Via a direct checking we get $d_{r}\equiv0$ for any $r\geq1$.
This implies that the Dolbeault--Koszul--Brylinski spectral sequence of
$(\mathbb{P}^{n}, \pi)$ degenerates at $E_{1}$-page, and therefore we obtain
$$
H_{k}(\mathbb{P}^{n},\pi)=
\left\{ \begin{array}{ll}
\mathbb{C}^{n+1},&k=n,\\
0, & k\neq n.
\end{array} \right.$$
Consider the Poisson blow-up $\varphi: (M_r,\tilde{\pi})\rightarrow(\mathbb{P}^{2},\pi)$.
From the blow-up formula in Theorem \ref{main-thm}, we get
$$
H_{k}(M_r,\tilde{\pi})=
\left\{ \begin{array}{ll}
\mathbb{C}^{r+3},&k=2,\\
0, & k\neq 2.
\end{array} \right.
$$

\subsection{Iwasawa manifolds}

To begin with, let us recall some basic facts on complex nilmanifolds.
Let $G$ be a complex nilpotent Lie group with Lie algebra $\mathfrak{g}$
whose complexification is $\mathfrak{g}_{\mathbb{C}}
:=\mathfrak{g}\mathfrak{}\otimes_{\mathbb{R}}\mathbb{C}$,
and let $H$ be a discrete subgroup of $G$.
Suppose $M=G/H$ is the associated nilmanifold endowed with a left-invariant
complex structure $J$ and a left-invariant holomorphic Poisson bi-vector field $\pi$.
Then there exists a natural inclusion of complexes
\begin{equation}\label{pos-in}
i: \bigl(\wedge^{p,\bullet}\mathfrak{g}_{\mathbb{C}}^{\ast},\bar{\partial}\bigr)
\hookrightarrow
\bigl(\Gamma(M,\A_M^{p,\bullet}),\bar{\partial}\bigr),
\end{equation}
for any $p\geq0$.
Set $n:=\dim_{\mathbb{C}}\, M$.

\begin{lem}\label{lem-inv-pos}
If the map \eqref{pos-in} is a quasi-isomorphism,
then the total cohomology of the double complex
$(\wedge^{\bullet,\bullet}\mathfrak{g}_{\mathbb{C}}^{\ast},\partial_\pi,\bar{\partial})$
is isomorphic to $H_{\bullet}(M,\pi)$.
\end{lem}
\begin{proof}
Observe that \eqref{pos-in} induces a morphism of double complexes
\begin{equation}\label{pos-doub}
i:
\bigl(\wedge^{\bullet,\bullet}\mathfrak{g}_{\mathbb{C}}^{\ast},\partial_\pi,\bar{\partial}\bigr)
\longrightarrow
\bigl(\Gamma(M,\A_M^{\bullet,\bullet}),\partial_\pi,\bar{\partial}\bigr).
\end{equation}
On the one hand, we know that
$(\wedge^{\bullet,\bullet}\mathfrak{g}_{\mathbb{C}}^{\ast},\partial_\pi,\bar{\partial})$
admits a spectral sequence $\{\tilde{E}_{r},\tilde{d}_{r}\}$
 converging to the corresponding total cohomology such that the $\tilde{E}_{1}$-page states as
\begin{equation*}
\tilde{E}_{1}^{p,q}= H^{q}(\wedge^{n-p,\bullet}\mathfrak{g}_{\mathbb{C}}^{\ast},\bar{\partial}).
\end{equation*}
On the other hand, the Dolbeault--Koszul--Brylinski  spectral sequence $\{E_{r},d_{r}\}$
converges to the holomorphic Koszul--Brylinski homology
$H_{\bullet}(M,{\pi})$ and has the $\tilde{E}_{1}$-page
\begin{equation*}
E_{1}^{p,q}=H^{q}(M,\Omega_M^{n-p}).
\end{equation*}
For any $r\geq1$, the inclusion \eqref{pos-doub} induces a morphism of the spectral sequences
$$
i_{r}^{\star}:\tilde{E}_{r}\longrightarrow E_{r}.
$$
Since \eqref{pos-in} is a quasi-isomorphism,
i.e., $i_1^{\star}:\tilde{E}_{1}\rightarrow E_{1}$ is an isomorphism,
by the standard result in the spectral sequence theory,
$i_r^{\star}$ is an isomorphism for any $r\geq2$.
This implies that the total cohomology of double complex
$(\wedge^{\bullet,\bullet}\mathfrak{g}_{\mathbb{C}}^{\ast},\partial_\pi,\bar{\partial})$
is isomorphic to the holomorphic Koszul--Brylinski homology $H_{\bullet}(M,{\pi})$.
\end{proof}

\begin{rem}
A result of Sakane \cite[Theorem 1]{Sak76} states that if a complex nilmanifold is complex parallelisable
(i.e., the holomorphic tangent bundle is holomorphically trivial),
then the inclusion \eqref{pos-in} is a quasi-isomorphism.
\end{rem}

Next we consider a concrete example.
Let $\mathrm{H}(3; \mathbb{C})$ be the Heisenberg Lie group:
$$
\mathrm{H}(3; \mathbb{C})=
\left\{ \small{\left(\begin{array}{ccc}
1 & z_{1} & z_{2}\\
0 & 1  & z_{3} \\
0 & 0 & 1
\end{array}
\right)}\Bigg| z_{1},z_2,z_3\in \mathbb{C} \right\}\subset \mathrm{GL}(3; \mathbb{C}).
$$
As a complex manifold, $\mathrm{H}(3; \mathbb{C})$ is isomorphic to $\mathbb{C}^{3}$.
Consider the discrete group
$\mathrm{G}_3:=\mathrm{Gl}(3; \mathbb{Z}[\sqrt{-1}])\cap \mathrm{H}(3; \mathbb{C})$,
where $\mathbb{Z}[\sqrt{-1}]=\{a+b\sqrt{-1}\mid a,b\in \mathbb{Z}\}$ is the Gaussian integers.
The left multiplication gives rise to a natural $\mathrm{G_3}$-action on $\mathrm{H}(3; \mathbb{C})$,
and the corresponding faithful $\mathrm{G}_3$-action on $\mathbb{C}^{3}$ is given by
$$
(a_{1}, a_{2}, a_{3})\cdot(z_{1}, z_{2}, z_{3}):=(z_{1}+a_{1}, z_{2}+a_{1}z_{3}+a_{2}, z_{3}+a_{3}),
$$
where $a_{1}, a_{2}, a_{3}\in \mathbb{Z}[\sqrt{-1}]$.
Such a $\mathrm{G}_3$-action yields a monomorphism $f: \mathrm{G}_3\to \mathrm{Aff}(\mathbb{C}^{3})$.
Here $\mathrm{Aff}(\mathbb{C}^{3})$ is the affine transformation group of $\mathbb{C}^{3}$.
Therefore, such a $\mathrm{G}_3$-action is properly discontinuous.
Furthermore, the $\mathrm{G}_3$-quotient space
$$
\mathbb{I}_3:=\mathbb{C}^{3}/ \mathrm{G}_3
$$
is a compact complex Calabi--Yau threefold,
called the {\it Iwasawa manifold}, which is non-K\"{a}hler, non-formal, and complex parallelisable.

Denote by $(\mathfrak{g}_{\mathbb{C}}^{\ast})^{1,0}$ the space of
left-invariant holomorphic differential forms on $\mathrm{H}(3; \mathbb{C})$.
Then $(\mathfrak{g}_{\mathbb{C}}^{\ast})^{1,0}$ has a basis:
$$
w^{1}=dz_{1},\,\,\, w^{2}=dz_{2}-z_{1}dz_{3},\,\,\, w^{3}=dz_{3},
$$
satisfying the structure equations:
$$
\begin{cases}
  w^{1}=0 ,   &   \\
  dw^{3}=0 ,   &   \\
  dw^{2}=-w^{1}\wedge w^{3}.    &
\end{cases}
$$
The  dual basis of Lie algebra of left-invariant holomorphic vector fields on
$\mathrm{H}(3; \mathbb{C})$, denoted by $\mathfrak{g}_{\mathbb{C}}^{1,0}$, is
$$
X_1=\frac{\partial}{\partial z_1},\,\,\,
X_2=\frac{\partial}{\partial z_2},\,\,\,
X_3=\frac{\partial}{\partial z_3}+z_1\frac{\partial}{\partial z_2}
$$
with the structure equations $[X_1,X_2]=[X_2,X_3]=0, [X_1,X_3]=X_2$.

Note that each left-invariant holomorphic bi-vector field $\pi$ on $\mathbb{I}_{3}$ is of the form
$\pi=c_1X_1\wedge X_2+c_2X_1\wedge X_3+c_3X_2\wedge X_3$, where $c_1,c_2$ and $c_3$ are constants.
In particular, a direct checking shows that $[\pi,\pi]=0$ holds if and only if $c_2=0$.
Since $\pi$ is left-invariant and $\mathbb{I}_{3}$ is  complex parallelisable,
by Lemma \ref{lem-inv-pos}, the holomorphic Koszul--Brylinski homology of
$(\mathbb{I}_{3},\pi)$ can be computed in terms of the total cohomology of the double complex
$(\wedge^{\bullet,\bullet}\mathfrak{g}^{\ast}_{\mathbb{C}},\partial_\pi,\bar{\partial})$.
Observe that $\pi$ is the linear combination of two compatible Poisson bi-vector fields
  $\pi_{12}=X_1\wedge X_2$ and $\pi_{23}=X_2\wedge X_3$.
Since $\partial_{\pi_{12}}=\partial_{\pi_{23}}=0$ we get $\partial_\pi=0.$
It follows that the Dolbeault--Koszul--Brylinski spectral sequence of $(\mathbb{I}_3,\pi)$
degenerates at $E_1$-page and therefore the Koszul--Brylinski homology
$H_{\bullet}(\mathbb{I}_{3},\pi)$ can be read off from the Hodge diamond of $\mathbb{I}_3$ (see figure below).
\begin{equation*}
\begin{array}{cccccc}
\begin{matrix}
&&&&&    1  \\
&&&&  3  &&  2 \\
&&& 3 && 6 && 2\\
&&1 && 6 &&  6 && 1\\
&&& 2 && 6 && 3\\
&&&&  2  &&  3 \\
&&&&&    1
\end{matrix}&\\
\;\;\;\;(\textrm{Hodge diamond of $\mathbb{I}_{3}$}) \\
\end{array}
\end{equation*}

As a result, we have the following table which records the holomorphic
Koszul--Brylinski homology of $(\mathbb{I}_{3},\pi)$.
 \begin{table}[!htbp]
 \centering
 \begin{tabular}{|p{2cm}|p{1cm}|p{1cm}|p{1cm}|p{1cm}|p{1cm}|p{1cm}|p{1cm}|}
\hline
\makecell[c]{$k$} & \makecell[c]{$0$} & \makecell[c]{$1$}  & \makecell[c]{$2$}
&\makecell[c]{$3$}
 & \makecell[c]{$4$}  & \makecell[c]{$5$}    &\makecell[c]{$6$}\\ \hline
\makecell[c]{$H_{k}(\mathbb{I}_3,\pi)$}
&  \makecell[c]{$\mathbb{C}$}  &  \makecell[c]{$\mathbb{C}^5 $}
 & \makecell[c]{$\mathbb{C}^{11} $ }   &\makecell[c]{$\mathbb{C}^{14}$}
 &  \makecell[c]{$\mathbb{C}^{11} $}     & \makecell[c]{$\mathbb{C}^{5} $ }
  &\makecell[c]{$\mathbb{C}$}             \\ \hline
\end{tabular}
\end{table}

Note that $(\mathbb{I}_3,\pi)$ is a unimodular holomorphic Poisson manifold.
Due to Proposition \ref{dual} and the isomorphism \eqref{Serre-Poincare-duality},
the holomorphic Koszul--Brylinski homology of $(\mathbb{I}_3,\pi)$ is
isomorphic to its Lichnerowicz--Poisson cohomology $H^{\bullet}(\mathbb{I}_{3},\pi)$.
So we obtain the following table.
\begin{table}[!htbp]
\centering
\begin{tabular}{|p{2cm}|p{1cm}|p{1cm}|p{1cm}|p{1cm}|p{1cm}|p{1cm}|p{1cm}|}
\hline
\makecell[c]{$k$} & \makecell[c]{$0$} & \makecell[c]{$1$}  & \makecell[c]{$2$}    &\makecell[c]{$3$}
 & \makecell[c]{$4$}  & \makecell[c]{$5$}    &\makecell[c]{$6$}\\ \hline
\makecell[c]{$H^{k}(\mathbb{I}_3,\pi)$}        &  \makecell[c]{$\mathbb{C}$}  &  \makecell[c]{$\mathbb{C}^5 $}     & \makecell[c]{$\mathbb{C}^{11} $ }   &\makecell[c]{$\mathbb{C}^{14}$}  &  \makecell[c]{$\mathbb{C}^{11} $}     & \makecell[c]{$\mathbb{C}^{5} $ }   &\makecell[c]{$\mathbb{C}$}             \\ \hline
\end{tabular}
\end{table}


\subsection{A six-dimensional complex nilmanifold}\label{a-six-nilfold}

Motivated by the construction of the Iwasawa manifold,
we consider the nilpotent Lie group
$$
G=\left\{A= \scriptsize{\left(\begin{array}{cccc}
1 &z_1 & z_{2} & z_{3}\\
0 &1 & z_{4} & z_{5}\\
0 & 0 & 1  & z_{6} \\
0& 0 & 0 & 1
\end{array}
\right)}\Bigg| z_{1}, z_2,\cdots,z_6\in \mathbb{C} \right\}\subset \mathrm{GL}(4; \mathbb{C}).
$$
As a complex manifold, $G$ is isomorphic to the complex vector space $\mathbb{C}^{6}$.
Consider the discrete subgroup $H:=\mathrm{Gl}(4; \mathbb{Z}[\sqrt{-1}])\cap G$.
Analogously, the left multiplication defines a natural $H$-action on $G$
and the corresponding faithful $H$-action on $\mathbb{C}^{6}$ is given by
\begin{multline*}
(a_{1}, a_{2}, a_{3},a_{4}, a_{5}, a_{6})\cdot(z_{1}, z_{2}, z_{3},z_4,z_5,z_6) \\
=(z_{1}+a_{1}, z_{2}+a_{1}z_{4}+a_{2},z_3+a_1z_5+a_2z_6+a_3,z_4+a_4,z_5+a_4z_6+a_5,z_6+a_6).
\end{multline*}
Therefore, this $H$-action is properly discontinuous, and the associated $H$-quotient space
$$
\mathbb{I}_6:=\mathbb{C}^{6}/ H
$$
is a compact complex manifold with complex dimension $6$.
Let $(\mathfrak{g}_{\mathbb{C}}^{\ast})^{1,0}$ be the space of
left-invariant holomorphic differential forms on $G$.
Then a basis of $(\mathfrak{g}_{\mathbb{C}}^{\ast})^{1,0}$ is given by
\begin{eqnarray*}
 &&w_1 =dz_1,\;\; w_2=dz_2-z_1dz_4,\;\; w_3=dz_3-z_1dz_5+(z_1z_4-z_2)dz_6,  \\
&& w_4 =dz_4,\;\; w_5=dz_5-z_4dz_6,\;\; w_6=dz_6. \qquad \qquad\qquad
\end{eqnarray*}
The structure equations are
$$
\begin{cases}
   dw_1=dw_4=dw_6=0,   &  \\
   dw_2=-w_1\wedge w_4,  &  \\
    dw_3=-w_1\wedge w_5-w_2\wedge w_6,     &  \\
     dw_5=-w_4\wedge w_6 .      &
\end{cases}
$$
Dually, Lie algebra of left-invariant holomorphic vector fields of $G$,
denoted by $\mathfrak{g}_{\mathbb{C}}^{1,0}$, has a basis:
\begin{eqnarray*}
&&X_{1}=\frac{\partial}{\partial z_1},
 \;\; X_2=\frac{\partial}{\partial z_2},
 \;\; X_3=\frac{\partial}{\partial z_3},\qquad \qquad\qquad \\
&&X_4=\frac{\partial}{\partial z_4}+z_1\frac{\partial}{\partial z_2},
 \;\; X_5=\frac{\partial}{\partial z_5}+z_1\frac{\partial}{\partial z_3},\;\;
 X_6=\frac{\partial}{\partial z_6}+z_2\frac{\partial}{\partial z_3}+z_4\frac{\partial}{\partial z_5}.
 \end{eqnarray*}
The only non-trivial relations of the dual basis are
$$
[X_1,X_4]=X_2,\,\,\, [X_1,X_5]=X_3=[X_2,X_6],\,\,\, [X_4,X_6]=X_5.
$$
It follows that $\mathbb{I}_6$ is a complex parallelisable,
non-K\"{a}hler, Calabi--Yau manifold with dimension $6$.

Now we consider some special holomorphic Poisson structures on
$\mathbb{I}_6$ given by left-invariant holomorphic bi-vector fields.
Akin to the Iwasawa manifold, the holomorphic Koszul--Brylinski homology of
$\mathbb{I}_6$ can be computed in terms of the total cohomology of the double complex
$(\wedge^{\bullet,\bullet}\mathfrak{g}^{\ast}_{\mathbb{C}},\partial_\pi,\bar{\partial})$.
For the simplicity, we write
$w^{i_{1}\cdots i_{p}\bar{j_{1}}\cdots \bar{j_{q}}}
=w^{i_{1}}\wedge \cdots \wedge w^{i_{p}}\wedge w^{\bar{j_{1}}}\wedge\cdots\wedge w^{\bar{j_{q}}}$,
for any $1\leq p,q\leq6$.
We study the holomorphic Koszul--Brylinski homology of $\mathbb{I}_6$ with respect to
the following three holomorphic Poisson bi-vector fields:
$$
\pi_{1}=X_2\wedge X_3,\;\;\;
\pi_{2}=X_1\wedge X_6,\;\;\;
\textrm{and} \;\;
\pi_{3}=X_1\wedge X_3.
$$

\subsubsection{Computation of $H_{\bullet}(\mathbb{I}_{6},\pi_{1})$}
We claim that the  the Dolbeault--Koszul--Brylinski spectral sequence of $(\mathbb{I}_{6}, \pi_{1})$
degenerates at $E_1$-page.
Observe that the only possible elements which are \emph{not}
$\partial_{\pi_1}$-closed are of the form
$w^{23i_{1}\cdots i_{p-2}\bar{j_{1}}\cdots \bar{j_{q}}}$.
However, a straightforward computation shows
\begin{eqnarray*}
     \partial_{\pi_1}w^{23i_{1}\cdots i_{p-2}\bar{j_{1}}\cdots \bar{j_{q}}}
&=&
     (\iota_{\pi_1}\circ \partial-\partial\circ \iota_{\pi_1})w^{23i_{1}\cdots
i_{p-2}\bar{j_{1}}\cdots \bar{j_{q}}}\\
     &=&\iota_{\pi_1}(w^{23}\wedge \partial w^{i_{1}\cdots i_{p-2}\bar{j_{1}}\cdots \bar{j_{q}}})
     -\partial w^{i_{1}\cdots i_{p-2}\bar{j_{1}}\cdots \bar{j_{q}}}\\
          &=&\partial w^{i_{1}\cdots i_{p-2}\bar{j_{1}}\cdots \bar{j_{q}}}
     -\partial w^{i_{1}\cdots i_{p-2}\bar{j_{1}}\cdots \bar{j_{q}}}\\
  &=&0.
  \end{eqnarray*}
This implies that the holomorphic volume form $\omega^{123456}$ is $\partial_{\pi_1}$-closed,
which means $(\mathbb{I}_6,\pi_1)$ is unimodular, and the
Dolbeault--Koszul--Brylinski spectral sequence of $(\mathbb{I}_6,\pi_1)$ degenerates at $E_1$-page.
Consequently, we get
\begin{equation}\label{degeneracy-I6}
H_k(\mathbb{I}_6,\pi_1)=\bigoplus\limits_{6-(p-q)=k} H_{\bar{\partial}}^{p,q}(\mathbb{I}_6).
\end{equation}
From the isomorphism \eqref{Serre-Poincare-duality}, we have
$$
H_k(\mathbb{I}_6,\pi_1)\cong H_{12-k}(\mathbb{I}_6,\pi_1).
$$
From the Hodge diamond of $\mathbb{I}_{6}$ (see Appendix \ref{appendix}) and Proposition \ref{dual},
we get the following table recording the holomorphic
Koszul--Brylinski homology of $(\mathbb{I}_6,\pi_1)$ up to degree 6
(the rest are obtained by the holomorphic Evens--Lu--Weinstein duality).
\begin{table}[!htbp]
\begin{tabular}{|p{2.3cm}|p{1cm}|p{1cm}|p{1cm}|p{1cm}|p{1cm}|p{1cm}|p{1cm}|}
\hline
\makecell[c]{$k$} & \makecell[c]{$0$} & \makecell[c]{$1$}  & \makecell[c]{$2$}    &\makecell[c]{$3$}
& \makecell[c]{$4$}  & \makecell[c]{$5$}    &\makecell[c]{$6$}\\ \hline
\makecell[c]{$H_{k}(\mathbb{I}_6,\pi_1)$}        &  \makecell[c]{$\mathbb{C}$}
&  \makecell[c]{$\mathbb{C}^9 $}     & \makecell[c]{$\mathbb{C}^{38} $ }
&\makecell[c]{$\mathbb{C}^{101}$}  &  \makecell[c]{$\mathbb{C}^{191} $}
& \makecell[c]{$\mathbb{C}^{274} $ }   &\makecell[c]{$\mathbb{C}^{308}$}             \\ \hline
\makecell[c]{$H^{12-k}(\mathbb{I}_6,\pi_1)$}        &  \makecell[c]{$\mathbb{C}$}
&  \makecell[c]{$\mathbb{C}^9 $}     & \makecell[c]{$\mathbb{C}^{38} $ }
&\makecell[c]{$\mathbb{C}^{101}$}  &  \makecell[c]{$\mathbb{C}^{191} $}
& \makecell[c]{$\mathbb{C}^{274} $ }   &\makecell[c]{$\mathbb{C}^{308}$}             \\ \hline
\end{tabular}
\end{table}


\begin{rem}
If the Dolbeault--Koszul--Brylinski spectral sequence for a holomorphic Poisson manifold
degenerates at the $E_{1}$-page,
then we can read off its holomorphic Koszul--Brylinski homology from the
Hodge diamond using the same method as in the
computation of $H_{\bullet}(\mathbb{I}_{6},\pi_{1})$.
However, the $E_{1}$-degeneracy of the Dolbeault--Koszul--Brylinski
spectral sequence is not a necessary condition
for a holomorphic Poisson manifold.
\end{rem}

\subsubsection{$E_{1}$-non-degeneracy for $(\mathbb{I}_{6},\pi_2)$}\label{ex1-non-degeeracy-E1}

Consider the holomorphic Poisson manifold $(\mathbb{I}_{6},\pi_2)$.
Observe that $(\mathfrak{g}_{\mathbb{C}}^{\ast})^{6,0}=\langle w^{123456}\rangle$ and
$\partial_{\pi_2} \omega^{123456}=0$; we obtain
\begin{equation*}\label{h0}
 H_{0}(\mathbb{I}_6, {\pi_2})=\langle [w^{123456}]\rangle\cong\mathbb{C}.
\end{equation*}
On the one hand, note that
$(\mathfrak{g}_{\mathbb{C}}^{\ast})^{5,0}=\langle w^{23456},w^{13456},w^{12456},w^{12356},w^{12346},w^{12345}\rangle$,
and we have
\begin{eqnarray*}
&&\partial_{\pi_2} w^{23456}=\partial_{\pi_2} w^{12456}=\partial_{\pi_2} w^{12345}=0, \\
   &&       \partial_{\pi_2} w^{13456}=-w^{2456},\\
    && \partial_{\pi_2} w^{12356}=w^{1345}-w^{2346},\\
    && \partial_{\pi_2} w^{12346}=-w^{1245}.
\end{eqnarray*}
On the other hand, since
$$(\mathfrak{g}_{\mathbb{C}}^{\ast})^{6,1}=\langle w^{123456\bar{1}},w^{123456\bar{2}},w^{123456\bar{3}},
                  w^{123456\bar{4}},w^{123456\bar{5}},w^{123456\bar{6}}\rangle,
$$
the following equalities hold:
$$\partial_{\pi_2}| (\mathfrak{g}_{\mathbb{C}}^{\ast})^{6,1}=0\,\,\,\,
\mathrm{and}
\,\,\,
\ker\,\bar{\partial}\cap(\mathfrak{g}_{\mathbb{C}}^{\ast})^{6,1}=
\langle w^{123456\bar{1}},
                  w^{123456\bar{4}},w^{123456\bar{6}}\rangle.
$$
Consequently, we get
\begin{center}
$H_{1}(\mathbb{I}_6,{\pi_2})=\langle [w^{23456}],[w^{12456}], [w^{12345}], [w^{123456\bar{1}}],
                  [w^{123456\bar{4}}],[w^{123456\bar{6}}]\rangle\cong\mathbb{C}^{6}.$
\end{center}
Assuming that the Dolbeault--Koszul--Brylinski spectral sequence of $(\mathbb{I}_6,\pi_2)$
degenerates at the $E_1$ page, we get
\begin{equation}\label{degeneracy-1}
H_1(\mathbb{I}_6,\pi_2)=H_{\bar{\partial}}^{5,0}(\mathbb{I}_6)\oplus H_{\bar{\partial}}^{6,1}(\mathbb{I}_6).
\end{equation}
Notice that $H_{\bar{\partial}}^{5,0}(\mathbb{I}_6)\cong \mathbb{C}^{6}$ and
$H_{\bar{\partial}}^{6,1}(\mathbb{I}_6)\cong \mathbb{C}^{3}$ (see Appendix \ref{appendix}).
This leads to a contradiction to the equality \eqref{degeneracy-1},
and therefore the Dolbeault--Koszul--Brylinski spectral sequence of $(\mathbb{I}_6,\pi_2)$
does not degenerate at the $E_1$-page.


\subsubsection{Computation of $H_{\bullet}(\mathbb{I}_{6},\pi_{3})$}

A direct computation shows that the non-trivial $\partial_{\pi_{3}}$-closed monomials are
given by:
\begin{itemize}
\item[(1)]On $(\mathfrak{g}_{\mathbb{C}}^{\ast})^{5,q}$, $\partial_{\pi_{3}}
w^{12356\bar{j_{1}}\cdots \bar{j_{q}}}
         =-w^{1456\bar{j_{1}}\cdots \bar{j_{q}}}$;
\item[(2)]On $(\mathfrak{g}_{\mathbb{C}}^{\ast})^{4,q}$,
\begin{eqnarray*}
   &&\partial_{\pi_{3}} w^{1235\bar{j_{1}}\cdots \bar{j_{q}}}=-w^{145\bar{j_{1}}\cdots \bar{j_{q}}} ,\\
   && \partial_{\pi_{3}} w^{1236\bar{j_{1}}\cdots \bar{j_{q}}}=-w^{146\bar{j_{1}}\cdots \bar{j_{q}}},\\
   && \partial_{\pi_{3}} w^{2356\bar{j_{1}}\cdots \bar{j_{q}}}=w^{456\bar{j_{1}}\cdots \bar{j_{q}}};
\end{eqnarray*}
\item[(3)]On $(\mathfrak{g}_{\mathbb{C}}^{\ast})^{3,q}$,
\begin{eqnarray*}
   &&\partial_{\pi_{3}} w^{123\bar{j_{1}}\cdots \bar{j_{q}}}=-w^{14\bar{j_{1}}\cdots \bar{j_{q}}},  \\
   &&\partial_{\pi_{3}} w^{235\bar{j_{1}}\cdots \bar{j_{q}}}=w^{45\bar{j_{1}}\cdots \bar{j_{q}}},\\
   && \partial_{\pi_{3}} w^{236\bar{j_{1}}\cdots \bar{j_{q}}}=w^{46\bar{j_{1}}\cdots \bar{j_{q}}};
\end{eqnarray*}
\item[(4)]On $(\mathfrak{g}_{\mathbb{C}}^{\ast})^{2,q}$,
$\partial_{\pi_{3}} w^{23\bar{j_{1}}\cdots \bar{j_{q}}}=w^{4\bar{j_{1}}\cdots \bar{j_{q}}}$.
\end{itemize}
It follows that $(\mathbb{I}_6,\pi_3)$ is unimodular, and the Dolbeault--Koszul--Brylinski
spectral sequence of $(\mathbb{I}_6,\pi_{3})$ does not
degenerate at the $E_1$-page.
By Lemma \ref{lem-inv-pos} and Proposition \ref{dual},
we have the following table.
\begin{table}[!htbp]
 \begin{tabular}{|p{2.3cm}|p{1cm}|p{1cm}|p{1cm}|p{1cm}|p{1cm}|p{1cm}|p{1cm}|}
\hline
\makecell[c]{$k$} & \makecell[c]{$0$} & \makecell[c]{$1$}  & \makecell[c]{$2$}
&\makecell[c]{$3$}
 & \makecell[c]{$4$}  & \makecell[c]{$5$}    &\makecell[c]{$6$}\\ \hline
\makecell[c]{$H_{k}(\mathbb{I}_6,\pi_3)$}        &  \makecell[c]{$\mathbb{C}$}
&  \makecell[c]{$\mathbb{C}^8 $}     & \makecell[c]{$\mathbb{C}^{31} $ }
&\makecell[c]{$\mathbb{C}^{78}$}  &  \makecell[c]{$\mathbb{C}^{143} $}
& \makecell[c]{$\mathbb{C}^{202} $ }   &\makecell[c]{$\mathbb{C}^{226}$}             \\ \hline
\makecell[c]{$H^{12-k}(\mathbb{I}_6,\pi_3)$}        &  \makecell[c]{$\mathbb{C}$}
&  \makecell[c]{$\mathbb{C}^8 $}     & \makecell[c]{$\mathbb{C}^{31} $ }
&\makecell[c]{$\mathbb{C}^{78}$}  &  \makecell[c]{$\mathbb{C}^{143} $}
& \makecell[c]{$\mathbb{C}^{202} $ }   &\makecell[c]{$\mathbb{C}^{226}$}             \\ \hline
\end{tabular}
\end{table}

\subsubsection{Poisson blow-up of $(\mathbb{I}_6,\pi_{3})$}

Take
$$\Gamma_2=
\left\{A= \scriptsize{\left(\begin{array}{cccc}
1 &z_{1} & z_{2} & z_{3}\\
0 &1 & a_{23} & a_{24}\\
0 & 0 & 1  & a_{34} \\
0& 0 & 0 & 1
\end{array}
\right)}\Bigg| z_{1},z_2,z_3\in \mathbb{C},a_{23},a_{24},a_{34}\in\mathbb{Z}[\sqrt{-1}] \right\}.
$$
Then $Y_2:=\Gamma_2/H$ is a $3$-dimensional K\"{a}hlerian nilmanfold.
Furthermore, $(Y_2,\pi_3|_{Y_2}=X_1\wedge X_3)$ is a closed holomorphic Poisson submanifold of
$(\mathbb{I}_6,\pi_{3})$ whose transverse Poisson structure vanishes.
One can check that $\partial_{\pi_3|_{Y_2}}=0$ and
thus the Dolbeault--Koszul--Brylinski spectral sequence of $(Y_{2},\pi_3|_{Y_2})$
degenerates at the $E_1$-page.
Note that the Hodge diamond of $Y_2$ is
\begin{equation*}
\begin{array}{cccccc}
\begin{matrix}
&&&&&    1  \\
&&&&  3  &&  3 \\
&&& 3 && 9 && 3\\
&&1 && 9 &&  9 && 1\\
&&& 3 && 9 && 3\\
&&&&  3  &&  3 \\
&&&&&    1
\end{matrix}&
\end{array}
\end{equation*}
As a corollary, we get the holomorphic Koszul--Brylinski homology of $(Y_{2},\pi_3|_{Y_2})$ as follows:
 \begin{table}[!htbp]
 \begin{tabular}{|p{2.5cm}|p{1cm}|p{1cm}|p{1cm}|p{1cm}|p{1cm}|p{1cm}|p{1cm}|}
\hline
\makecell[c]{$k$} & \makecell[c]{$0$} & \makecell[c]{$1$}  & \makecell[c]{$2$}
&\makecell[c]{$3$}
 & \makecell[c]{$4$}  & \makecell[c]{$5$}    &\makecell[c]{$6$}\\ \hline
\makecell[c]{$H_{k}(Y_2, \pi_3|_{Y_2})$}        &  \makecell[c]{$\mathbb{C}$}
&  \makecell[c]{$\mathbb{C}^6 $}     & \makecell[c]{$\mathbb{C}^{15} $ }
&\makecell[c]{$\mathbb{C}^{20}$}  &  \makecell[c]{$\mathbb{C}^{15} $}
& \makecell[c]{$\mathbb{C}^{6} $ }   &\makecell[c]{$\mathbb{C}^{1}$}              \\ \hline
\end{tabular}
\end{table}

Let $\varphi:\mathrm{Bl}_{Y_2}\mathbb{I}_6\rightarrow\mathbb{I}_{6}$ be the blow-up of $\mathbb{I}_{6}$
along $Y_{2}$.
By Proposition \ref{Constr-Poisson-bp}, the holomorphic Poisson structure $\pi_{3}$ can be lifted to
a unique holomorphic Poisson structure $\tilde{\pi}_{3}$ on $\mathrm{Bl}_{Y_2}\mathbb{I}_6$.
By Theorem \ref{main-thm}, we get the  the following table recording
the holomorphic Koszul--Brylinski homology of $(\mathrm{Bl}_{Y_2}\mathbb{I}_6,\tilde{\pi}_3)$.
 \begin{table}[!htbp]
 \begin{tabular}{|p{2.5cm}|p{1cm}|p{1cm}|p{1cm}|p{1cm}|p{1cm}|p{1cm}|p{1cm}|}
\hline
\makecell[c]{$k$} & \makecell[c]{$0$} & \makecell[c]{$1$}  & \makecell[c]{$2$}    &\makecell[c]{$3$}
 & \makecell[c]{$4$}  & \makecell[c]{$5$}    &\makecell[c]{$6$}\\ \hline
\makecell[c]{$H_{k}(\mathrm{Bl}_{Y_2}\mathbb{I}_6,\tilde{\pi}_3)$}
&  \makecell[c]{$\mathbb{C}$}  &  \makecell[c]{$\mathbb{C}^8$}
& \makecell[c]{$\mathbb{C}^{31} $ }   &\makecell[c]{$\mathbb{C}^{80}$}
&  \makecell[c]{$\mathbb{C}^{155} $}     & \makecell[c]{$\mathbb{C}^{232} $ }
&\makecell[c]{$\mathbb{C}^{266}$}             \\ \hline
\end{tabular}
\end{table}


\appendix
\section{Hodge diamond of $\mathbb{I}_6$}\label{app-I6}\label{appendix}

Note that $\mathbb{I}_6$ is complex parallelisable.
As mentioned in the main text,  the Dolbeault cohomology of $\mathbb{I}_6$
can be computed by means of left-invariant forms (\cite[Theorem 1]{Sak76}).
Consider the associated double complex
$(\wedge^{\bullet,\bullet}\mathfrak{g}_{\mathbb{C}}^{\ast},\partial,\bar{\partial})$.
By Leibniz rule, we have
$$
\bar{\partial}w^{i_{1}\cdots i_{p}\bar{j_1}\cdots\bar{j_q}}
=(-1)^{p}w^{i_{1}\cdots i_{p}}\wedge\bar{\partial}w^{\bar{j_1}\cdots\bar{j_q}}.
$$
In particular, we get  $h^{i,j}=\tbinom{6}{i}\cdot h^{0,j}$,
where $h^{i,j}:=\mathrm{dim}_{\mathbb{C}}\,
H^{j}(\wedge^{i,\bullet}\mathfrak{g}_{\mathbb{C}}^{\ast},\bar{\partial})$
is the Lie algebra Hodge number.
For this reason, to compute the Hodge diamond of $\mathbb{I}_6$,
we only need to compute $h^{0,0}, h^{0,1},\cdots, h^{0,6}$.
Since $\mathbb{I}_6$ is a compact complex manifold we have
$H_{\bar{\partial}}^{0,0}(\mathbb{I}_6)\cong\mathbb{C}.$
The monomials in $(\mathfrak{g}_{\mathbb{C}}^{\ast})^{0,j}$
which are \emph{not} $\bar{\partial}$-closed are stated as follows:
\begin{itemize}
\item[(1)] On $(\mathfrak{g}_{\mathbb{C}}^{\ast})^{0,1}$, since
$$\bar{\partial}w^{\bar{2}}=-w^{\bar{1}\bar{4}}, \quad
    \bar{\partial}w^{\bar{3}}=-w^{\bar{1}\bar{5}}-w^{\bar{2}\bar{6}},\quad
     \bar{\partial}w^{\bar{5}}=-w^{\bar{4}\bar{6}},
$$
we get
\begin{center}
$H_{\bar{\partial}}^{0,1}(\mathbb{I}_6)
      =\langle [w^{\bar{1}}],[w^{\bar{4}}], [w^{\bar{6}}]\rangle\cong\mathbb{C}^{3}.$
\end{center}

\item[(2)] On $(\mathfrak{g}_{\mathbb{C}}^{\ast})^{0,2}$, since
\begin{eqnarray*}
&&\bar{\partial}w^{\bar{1}\bar{3}}=w^{\bar{1}\bar{2}\bar{6}}, \quad
\bar{\partial}w^{\bar{1}\bar{5}}=w^{\bar{1}\bar{4}\bar{6}},  \quad
\bar{\partial}w^{\bar{2}\bar{3}}=w^{\bar{1}\bar{3}\bar{4}}-w^{\bar{1}\bar{2}\bar{5}},\\
&&\bar{\partial}w^{\bar{2}\bar{6}}=-w^{\bar{1}\bar{4}\bar{6}},\quad
\bar{\partial}w^{\bar{3}\bar{6}}=-w^{\bar{1}\bar{5}\bar{6}},\quad
\bar{\partial}w^{\bar{3}\bar{4}}=w^{\bar{1}\bar{4}\bar{5}} +w^{\bar{2}\bar{4}\bar{6}},\\
&&\bar{\partial}w^{\bar{2}\bar{5}}=-w^{\bar{1}\bar{4}\bar{5}}+w^{\bar{2}\bar{4}\bar{6}},
\quad \bar{\partial}w^{\bar{3}\bar{5}}=w^{\bar{2}\bar{5}\bar{6}}+w^{\bar{3}\bar{4}\bar{6}},
  \end{eqnarray*}
we get
\begin{center}
$H_{\bar{\partial}}^{0,2}(\mathbb{I}_6)
=\langle [w^{\bar{1}\bar{2}}],[w^{\bar{1}\bar{6}}],[w^{\bar{2}\bar{4}}],
[w^{\bar{4}\bar{5}}],[w^{\bar{5}\bar{6}}]\rangle\cong\mathbb{C}^{5}.$
\end{center}

\item[(3)] On $(\mathfrak{g}_{\mathbb{C}}^{\ast})^{0,3}$, since
\begin{eqnarray*}
&&\bar{\partial}w^{\bar{1}\bar{2}\bar{5}}=-w^{\bar{1}\bar{2}\bar{4}\bar{6}},\quad
\bar{\partial}w^{\bar{1}\bar{3}\bar{4}}=-w^{\bar{1}\bar{2}\bar{4}\bar{6}},  \quad
\bar{\partial}w^{\bar{2}\bar{3}\bar{4}} =w^{\bar{1}\bar{2}\bar{4}\bar{5}},\\
&&\bar{\partial}w^{\bar{2}\bar{5}\bar{6}}=-w^{\bar{1}\bar{4}\bar{5}\bar{6}},\quad
\bar{\partial}w^{\bar{3}\bar{4}\bar{5}}=-w^{\bar{2}\bar{4}\bar{5}\bar{6}},\quad
\bar{\partial}w^{\bar{3}\bar{4}\bar{6}}=w^{\bar{1}\bar{4}\bar{5}\bar{6}},\\
&&\bar{\partial}w^{\bar{1}\bar{3}\bar{5}}=-w^{\bar{1}\bar{2}\bar{5}\bar{6}}
             -w^{\bar{1}\bar{3}\bar{4}\bar{6}},\quad
\bar{\partial}w^{\bar{2}\bar{3}\bar{5}}
       =w^{\bar{1}\bar{3}\bar{4}\bar{5}} -w^{\bar{2}\bar{3}\bar{4}\bar{6}},\quad
\bar{\partial}w^{\bar{2}\bar{3}\bar{6}}=w^{\bar{1}\bar{3}\bar{4}\bar{6}}
                  -w^{\bar{1}\bar{2}\bar{5}\bar{6}},
  \end{eqnarray*}
we get
\begin{center}
$H_{\bar{\partial}}^{0,3}(\mathbb{I}_6)
=\langle [w^{\bar{1}\bar{2}\bar{4}}],[w^{\bar{1}\bar{3}\bar{4}}],[w^{\bar{1}\bar{3}\bar{6}}],
[w^{\bar{2}\bar{4}\bar{5}}],[w^{\bar{3}\bar{5}\bar{6}}],[w^{\bar{4}\bar{5}\bar{6}}]\rangle
\cong\mathbb{C}^{6}.$
\end{center}

\item[(4)] On $(\mathfrak{g}_{\mathbb{C}}^{\ast})^{0,4}$, since
\begin{eqnarray*}
&&\bar{\partial}w^{\bar{1}\bar{2}\bar{3}\bar{5}}
      =w^{\bar{1}\bar{2}\bar{3}\bar{4}\bar{6}},  \quad
\bar{\partial}w^{\bar{1}\bar{3}\bar{4}\bar{5}}=w^{i_{1}\cdots i_{p}\bar{1}\bar{2}\bar{4}\bar{5}\bar{6}},  \quad
\bar{\partial}w^{\bar{2}\bar{3}\bar{4}\bar{6}}=w^{\bar{1}\bar{2}\bar{4}\bar{5}\bar{6}},\quad
\bar{\partial}w^{\bar{2}\bar{3}\bar{5}\bar{6}}=w^{\bar{1}\bar{3}\bar{4}\bar{5}\bar{6}},
  \end{eqnarray*}
we get
\begin{center}
$H_{\bar{\partial}}^{0,4}(\mathbb{I}_6)
=\langle [w^{\bar{1}\bar{2}\bar{3}\bar{4}}],[w^{\bar{1}\bar{2}\bar{3}\bar{6}}],
[w^{\bar{1}\bar{3}\bar{5}\bar{6}}],[w^{\bar{2}\bar{3}\bar{4}\bar{5}}],[w^{\bar{3}\bar{4}\bar{5}\bar{6}}]
\rangle\cong\mathbb{C}^{5}.$
\end{center}
\item[(5)]Since $\bar{\partial}|_{(\mathfrak{g}_{\mathbb{C}}^{\ast})^{0,5}}
=\bar{\partial}|_{(\mathfrak{g}_{\mathbb{C}}^{\ast})^{0,6}}=0$,
we have
$$H_{\bar{\partial}}^{0,5}(\mathbb{I}_6)=
\langle [w^{\bar{1}\bar{2}\bar{3}\bar{4}\bar{5}}],[w^{\bar{1}\bar{2}\bar{3}\bar{5}\bar{6}}]
[w^{\bar{2}\bar{3}\bar{4}\bar{5}\bar{6}}]\rangle\cong\mathbb{C}^3\quad
\textup{and}\quad
H_{\bar{\partial}}^{0,6}(\mathbb{I}_6)=
\langle [w^{\bar{1}\bar{2}\bar{3}\bar{4}\bar{5}\bar{6}}]\rangle\cong\mathbb{C}.$$
\end{itemize}
By the discussion in the above,
we obtain the Hodge diamond of $\mathbb{I}_6$ as follows:
\begin{equation*}
\small{
   \begin{array}{ccccccccccccc}
       &   &   &   &   &   & 1  &   &   &   &   &   &  \\
       &   &   &   &   & 6 &    & 3 &   &   &   &   &  \\
       &   &   &   & 15&   & 18 &   & 5 &   &   &   &  \\
       &   &   &20 &   &45 &    &30 &   & 6 &   &   &  \\
       &   &15 &   &60 &   &75  &   &36 &   & 5 &   &  \\
       & 6 &   &45 &   &100&    &90 &   &30 &   & 3 &  \\
     1 &   &18 &   &75 &   &120 &   &75 &   &18 &   & 1 \\
       & 3 &   &30 &   &90 &    &100&   &45 &   & 6 &  \\
       &   & 5 &   &36 &   &75  &   &60 &   &15 &   &  \\
       &   &   & 6 &   &30 &    &45 &   &20 &   &   &  \\
       &   &   &   & 5 &   &18  &   &15 &   &   &   &  \\
       &   &   &   &   & 3 &    & 6 &   &   &   &   &  \\
       &   &   &   &   &   & 1  &   &   &   &   &   &  \\
   \end{array}}
\end{equation*}



\begin{thebibliography}{GP}

\bibitem{ASTT20}
D.\ Angella, T.\ Suwa, N.\ Tardini, and A.\ Tomassini,
{\it Note on Dolbeault cohomology and Hodge structures up to bimeromorphisms},
Complex Manifolds  7 (2020) 194--214.

\bibitem{Ba13}
M. Bailey,
{\it Local classification of generalized complex structures},
J. Differential Geom. 95 (2013), 1--37.


\bibitem{BCV19}
M. Bailey, G. R. Cavalcanti, and J. L. van der Leer Dur\'{a}n,
{\it Blow-ups in generalized complex geometry},
Trans. Amer. Math. Soc.  371 (2019),  2109--2131.


\bibitem{Bon93}
A. Bondal,
{\it Non-commutative deformations and Poisson brackets on projective spaces},
Preprint no. 67, Max-Planck-Institut, Bonn 1993.


\bibitem{BX15}
D. Broka and P. Xu,
{\it Symplectic realizations of holomorphic Poisson manifolds},
to appear in Math. Res. Lett. arXiv:1512.08847.


\bibitem{Bry88}
J.-L. Brylinski,
{\it A differential complex for Poisson manifolds},
J. Differential Geom.  28 (1988), 93--114.


\bibitem{BZ99}
J.-L. Brylinski and G. Zuckerman,
{\it The outer derivation of a complex Poisson manifold},
J. Reine Angew. Math. 506 (1999), 181--189.


\bibitem{CFP16}
Z. Chen, A. Fino, and Y.-S. Poon,
{\it Holomorphic Poisson structure and its cohomology on nilmanifolds},
Differential Geom. Appl. 44 (2016), 144--160.


\bibitem{CSX10}
Z. Chen, M. Sti\'{e}non, and P. Xu,
{\it Geometry of Maurer--Cartan elements on complex Manifolds},
Comm. Math. Phys. 297 (2010),  169--187.


\bibitem{CY22}
Y. Chen and S. Yang,
{\it On blow-up formula of integral Bott--Chern cohomology},
Ann. Glob. Anal. Geom. 61 (2022), 57--67.

\bibitem{CGP15}
Z. Chen, D. Grandini, and Y.-S. Poon,
{\it Holomorphic Poisson cohomology},
Complex Manifolds 2 (2015), 34--52.


\bibitem{DGMS75}
P. Deligne, P. Griffiths, J. Morgan, and D. Sullivan,
{\it Real homotopy theory of K\"{a}hler manifolds},
Invent. Math. 29 (1975), 245--274.



\bibitem{Dem12}
J.-P. Demailly,
{\it Complex analytic and differential geometry}, available at
\href{https://www-fourier.ujf-grenoble.fr/~demailly/manuscripts/agbook.pdf}
{https://www-fourier.ujf-grenoble.fr/$\sim$demailly/manuscripts/agbook.pdf}.


\bibitem{FM12}
D. Fiorenza and M. Manetti,
{\it Formality of Koszul brackets and deformations of holomorphic Poisson manifolds},
Homology, Homotopy Appl. 14 (2012), 63--75.


\bibitem{Fu05}
B. Fu,
{\it Poisson resolutions},
J. reine angew. Math. 587  (2005), 17--26.


\bibitem{Go10}
R. Goto,
{\it Deformations of generalized complex and generalized K\"{a}hler structures},
J. Differential Geom. 84 (2010), 525--560.



\bibitem{Gu11}
M. Gualtieri,
{\it Generalized complex geometry},
 Ann. of Math.  174 (2011), 75--123.

\bibitem{Hit03}
N.J. Hitchin,
{\it Generalized Calabi-Yau manifolds},
Quart. J. Math. 54 (2003), 281--308.


\bibitem{Hit06}
N.J. Hitchin,
{\it Instantons, Poisson structures and generalized K\"{a}hler geometry},
Comm. Math. Phys. 265 (2006), 131--164.


\bibitem{Hit12}
N.J. Hitchin,
{\it Deformations of holomorphic Poisson manifolds},
Mosc. Math. J. 669 (2012), 567--591.


\bibitem{Hon19}
W. Hong,
{\it Poisson cohomology of holomorphic toric Poisson manifolds. I},
J. Algebra  527 (2019), 147--181.


\bibitem{HX11}
W. Hong and P. Xu,
{\it Poisson cohomology of Del Pezzo surfaces},
J. Algebra 336 (2011), 378--390.


\bibitem{Kod86}
K. Kodaira,
{\it Complex manifolds and deformation of complex structures},
Classics in Mathematics, Springer, Berlin (2005).


\bibitem{Kos84}
J. L. Koszul,
{\it Crochet de Schouten-Nijenhuis et cohomologie},
The mathematical heritage of \'{E}lie Cartan (Lyon, 1984).


\bibitem{L-GSX08}
C. Laurent-Gengoux, M. Sti\'{e}non, and P. Xu,
{\it Holomorphic Poisson manifolds and holomorphic Lie algebroids},
Int.  Math. Res. Not. 2008, Art. ID rnn 088, 46 pp.


\bibitem{L-GPV13}
C. Laurent-Gengoux, A. Pichereau, and P. Vanhaecke,
{\it Poisson Structures},
Grundlehren Math. Wiss. 347, Springer, Heidelberg, (2013).


\bibitem{Lic77}
A. Lichnerowicz,
{\it Les vari\'{e}t\'{e}s de Poisson et leurs alg\`{e}bres de Lie associ\'{e}es},
J. Differential Geom. 12  (1977), 253--300.

\bibitem{Mc01}
J.\ McCleary,
{\it A user's guide to spectral sequences},
Second edition. Cambridge Studies in Advanced Mathematics, 58.
Cambridge University Press, Cambridge, (2001).

\bibitem{Pol97}
A. Polishchuk,
{\it Algebraic geometry of Poisson brackets},
Algebraic geometry, 7. J. Math. Sci. (New York),  84 (1997), 1413--1444.


\bibitem{PS17}
Y.-S. Poon and J. Simanyi,
{\it A Hodge-type decomposition of holomorphic Poisson cohomology on nilmanifolds},
Complex Manifolds 4 (2017), 137--154.


\bibitem{PS19}
Y.-S. Poon and J. Simanyi,
{\it Algebraic structure of holomorphic Poisson cohomology on nilmanifolds},
Complex Manifolds 6 (2019), 88--102.

\bibitem{Pym18}
B. Pym,
{\it Constructions and classifications of projective Poisson varieties},
Lett. Math. Phys. 108 (2018), 573--632.


\bibitem{RYY19}
S. Rao, S. Yang, and X. Yang,
{\it Dolbeault cohomologies of blowing up complex manifolds},
J. Math. Pures Appl. 130 (2019), 68--92.


\bibitem{RYY20}
S. Rao, S. Yang, and X. Yang,
{\it Dolbeault cohomologies of blowing up complex manifolds II: bundle-valued case},
J. Math. Pures Appl. 133 (2020), 1--38.



\bibitem{Sak76}
Y. Sakane,
{\it On compact complex parallelisable solvmanifolds},
Osaka J. Math. 13 (1976), 187--212.

\bibitem{ST08}
G.\, Sharygin and D.\, Talalaev,
{\it On the Lie-formality of Poisson manifolds},
J. K-Theory 2 (2008) 361--384.

\bibitem{Sti11}
M. Sti\'{e}non,
{\it Holomorphic Koszul--Brylinski Homology},
Int.  Math. Res. Not.  2011 (2011), 553--571.


\bibitem{We97}
A. Weinstein,
{\it The modular automorphism group of a Poisson manifold},
J. Geom.  Phys.  23 (1997), 379--394.


\bibitem{YY20}
S. Yang and X. Yang,
{\it Bott--Chern blow-up formulae and the bimeromorphic invariance of the
$\partial\bar{\partial}$-lemma for threefolds},
Trans. Amer. Math. Soc.  373 (2020), 8885--8909.

\end{thebibliography}
\end{document}